\documentclass[10pt,reqno]{amsart}

\usepackage{bbm}
\usepackage{amsmath}
\usepackage{a4wide}

\font\cmssl=cmss10 at 10 pt  

\newtheorem{theorem}{Theorem}[section]
\newtheorem{lemma}[theorem]{Lemma}
\newtheorem{Prop}[theorem]{Proposition}
\newtheorem{Cor}[theorem]{Corollary}
\theoremstyle{definition}
\newtheorem{definition}[theorem]{Definition}
\newtheorem{example}[theorem]{Example}

\theoremstyle{remark}
\newtheorem{remark}[theorem]{Remark}

\numberwithin{equation}{section}

\newcommand{\id}   {{\mathbbm{1}}}

\newcommand{\bt}{\begin{theorem}\ \ }  
\newcommand{\et}{\end{theorem}}  
\newcommand{\bp}{\begin{Prop}\ \ }  
\newcommand{\ep}{\end{Prop}}  
\newcommand{\bc}{\begin{Cor}\ \ }  
\newcommand{\ec}{\end{Cor}}  
\newcommand{\bl}{\begin{lemma}\ \ }  
\newcommand{\el}{\end{lemma}}  
\newcommand{\bd}{\begin{definition}\ \ }  
\newcommand{\ed}{\end{definition}}  
\newcommand{\pf}{\begin{proof}}  
\newcommand{\epf}{\end{proof}}  
\newcommand{\br}{\begin{remark}\ \ }
\newcommand{\er}{\end{remark}}
\newcommand{\brsn}{\begin{remarks*}\ \ }
\newcommand{\ersn}{\end{remarks*}}

\newcommand{\be}{\begin{equation}}
\newcommand{\ee}{\end{equation}}

\newcommand{\arr}{\begin{array}{rlll}}
\newcommand{\ea}{\end{array}}
\newcommand{\bea}{\begin{eqnarray}}
\newcommand{\eea}{\end{eqnarray}}
\newcommand{\bean}{\begin{eqnarray*}}
\newcommand{\eean}{\end{eqnarray*}}

\newcommand{\ra}{\rightarrow}

\newcommand{\n}{\nabla}

\newcommand{\bR}{\mathbb{R}}
\newcommand{\bC}{\mathbb{C}}
\renewcommand{\o}{\omega}

\begin{document}
\title{On the structure of nearly pseudo-K\"ahler manifolds}

\author{Lars Sch\"afer}
\address{Lars Sch\"afer, Institut Differentialgeometrie, Leibniz Universit\"at
  Hannover,  Welfengarten 1, D-30167 Hannover, Germany}
\email{schaefer@math.uni-hannover.de}
\date{December 11, 2009}
\begin{abstract}
Firstly we give a condition to split off the K\"ahler factor from a nearly
pseudo-K\"ahler manifold and apply this to get a structure result in dimension
$8.$ Secondly we extend the construction of nearly K\"ahler manifolds from twistor spaces
to negatively curved quaternionic K\"ahler manifolds and para-quaternionic K\"ahler
manifolds. The class of nearly pseudo-K\"ahler manifolds obtained from this construction is
characterized by a holonomic condition. 
The combination of these results enables us to give a classification  result in (real) dimension 10. 
Moreover, we show that a strict nearly pseudo-K\"ahler six-manifold is Einstein. 
\end{abstract}

\maketitle
\section*{Introduction}
Nearly K\"ahler geometry was introduced and studied in a series of papers by A.\ Gray in the seventies in the context of weak holonomy. To our
best knowledge he only considers pseudo-Riemannian metrics in his paper on 3-symmetric spaces \cite{G2}. 
In the analysis of Killing spinors on pseudo-Riemannian manifolds \cite{Ka}
nearly pseudo-K\"ahler and nearly para-K\"ahler manifolds appear in a natural
way.  Levi-Civita flat nearly K\"ahler manifolds provide  a special class of solutions of the topological-antitopological
fusion equations on the tangent bundle \cite{S1,S2}. There is a strong similarity to
special K\"ahler geometry. For these reasons we became interested in Levi-Civita flat
nearly K\"ahler manifolds and were able to give a constructive classification
\cite{CS,CS2}.  In particular it follows, that non-K\"ahlerian examples only
exist in pseudo-Riemannian geometry and that the real dimension is at least 12. In other words, nearly K\"ahler geometry in the pseudo-Riemannian
world can be very different from the better-understood Riemannian world. In
(real) dimension six nearly pseudo-K\"ahler manifolds satisfy an exterior system analogue to the Riemannian case. Details can be found in
 \cite{SSH}. This system is used there to study such structures on products $G\times G,$ where $G$ is a simple three-dimensional Lie group. 

An interesting class of nearly K\"ahler manifolds $M^{4n+2}$ can be constructed from twistor
spaces over positive quaternionic K\"ahler manifolds. This class is
characterized \cite{BM,N1} by the reducibility  of the holonomy of the canonical
connection $\bar \n$ to $U(n) \times U(1)$. We show in this article that the twistor
spaces over negative quaternionic K\"ahler manifolds and para-quaternionic
K\"ahler manifolds carry a nearly pseudo-K\"ahler structure and characterize
the class of such examples by a holonomic condition. 

Using this result, we classify {\it nice and decomposable} (cf. Definitions
\ref{nice_def} and \ref{decomp_def}) nearly pseudo-K\"ahler manifolds in
dimension ten. \\

\noindent 
{\bf Theorem A.}
Let $(M^{10},J,g)$ be a nice decomposable nearly K\"ahler manifold,
then the universal cover of $M$ is either the product of a pseudo-K\"ahler 
surface and a (strict) nearly pseudo-K\"ahler manifold $M^6$ or a twistor
space of an eight-dimensional (para-)quaternionic K\"ahler manifold endowed 
with its canonical nearly pseudo-K\"ahler structure. \\

A strict nearly pseudo-K\"ahler six-manifold $M^6$ is shown to be Einstein in Theorem \ref{Einst_dim_6}. In dimension eight simply connected strict nearly
pseudo-K\"ahler manifolds are shown (Theorem \ref{class_dim8}) to be
products $\Sigma \times M^6$ of a Riemannian surface $\Sigma$ and a strict
nearly pseudo-K\"ahler manifold $M^6.$ 

In the first section of this paper we recall the definition of a nearly pseudo-K\"ahler manifold and generalize 
some facts and curvature identities to arbitrary signature. In the second
section we give a general condition to split off the K\"ahler factor from a
nearly pseudo-K\"ahler manifold, see Theorem \ref{Kaehler_factor_prop}. Using some linear algebra of three-forms this 
shows the splitting result for {\it nice} nearly pseudo-K\"ahler manifolds in 
dimension $8.$ The argument also holds true for a Riemannian metric and 
gives an alternative proof for the known result. If a {\it nice} nearly 
pseudo-K\"ahler ten-manifold is in addition {\it decomposable}, we find two 
cases: In the first we can split off the K\"ahler factor and in the second 
the holonomy of $\bar \n$ is reducible with a complex one dimensional factor. 
This is one motivation to study twistor spaces. Before doing this in section
four  we recall some information on pseudo-Riemannian submersions in the third section.
 In the pseudo-Riemannian setting twistor spaces are a good source of
 examples, since quaternionic geometry is richer in negative scalar curvature
 than in positive (cf. Remark \ref{REM_Example_twist}) and since we have the additional class of twistor spaces over para-quaternionic manifolds. In the last section we prove
that a nearly pseudo-K\"ahler manifold $M$ of {\it twistorial type}
(cf. Definition \ref{def_twist_type}) is obtained from the
above mentioned construction on a twistor space. This is done as follows: We prove that
$M$ comes from a pseudo-Riemannian submersion $\pi \,:\, M \ra N.$   Then we use
the nearly K\"ahler data on $M$ to endow $N$ with the structure of a (para-)quaternionic manifold. 
The proof is finished by identifying the twistor space of $N$ with $M.$ The
former proofs \cite{BM,N1} in the Riemannian case all use the inverse twistor construction
of Penrose or LeBrun, which does not seem to be developed for the situations
occurring in this text. As the reader might observe, the approach presented
here holds also true for Riemannian metrics. 

\noindent 
{\sl Acknowledgments.}
The author thanks Vicente Cort\'es for discussions.

\section{Nearly pseudo-K\"ahler manifolds}
\bd
An almost pseudo-hermitian manifold $(M,J,g)$ is called {\cmssl nearly pseudo-K\"ahler manifold} if it holds
$$(\n_XJ)X =0, \quad \forall X \in \Gamma(TM),$$
where $\n$ is the Levi-Civita connection of the (pseudo-)Riemannian metric
$g.$ A nearly pseudo-K\"ahler manifold is called {\cmssl strict} 
if it holds $\n_XJ\ne 0$ for all $X \in TM.$
\ed

\subsection{Curvature identities in the pseudo-Riemannian case}
The starting point of a series of curvature identities are 
\bea 
R(W,X,Y,Z)&-&R(W,X,JY,JZ) =g((\n_W J)X, (\n_Y J)Z ),
 \label{curv_id_i_equ} \\
R(W,X,W,Z)&+& R(W,JX,W,JZ) \label{curv_id_ii_equ} \\
&-&R(W,JW,X,JZ)= 2 g((\n_WJ)X, (\n_W J)Z ), \nonumber\\
R(W,X,Y,Z)&=&R(JW,JX,JY,JZ), \label{type_cond}
\eea
which were already proven  for pseudo-Riemannian metrics by Gray \cite{G1}.
Let $\{e_i\}_{i=1}^{2n}$ be a local orthonormal frame field, then
the Ricci- and the Ricci*-tensor are given by
\bean
g(Ric\, X,Y)= \sum_{i=1}^{2n} \epsilon_i\, R(X,e_i,Y,e_i), \quad 
g(Ric^*\, X,Y)= \frac{1}{2} \sum_{i=1}^{2n} \epsilon_i \, R(X,JY,e_i,Je_i)
\eean
with $\epsilon_i=g(e_i,e_i)=g(Je_i,Je_i)$ and $X,Y \in TM.$ The frame
$\{e_i\}_{i=1}^{2n}$ is called {\cmssl adapted} if it holds $Je_i=e_{i+n}$ for
$i=1,\ldots, n.$ Then it follows using an adapted frame
from equations \eqref{curv_id_ii_equ} and \eqref{type_cond} that 
\bea g(rX,Y):= g((Ric-Ric^*)X,Y)=   \sum_{i=1}^{2n} \epsilon_i\, 
g((\n_XJ)e_i, (\n_Y J)e_i ).
\eea  
Using the right hand-side we see $$[J, r]=0.$$ 
For the second derivative of the complex structure one has the identity
\bea \label{id_equ_D2_J} 2g(\n^2_{W,X}(J)Y,Z)= -\sigma_{X,Y,Z}\, g((\n_WJ)X,(\n_YJ)JZ), \eea
which was proven in \cite{G1} for Riemannian metrics and holds true in the pseudo-Riemannian setting, cf.  \cite{Ka} Proposition
7.1. This identity implies
\bea \label{id_equ_spur_D2_J} &~&\sum_{i=1}^{2n} \epsilon_i \n_{e_i,e_i}^2(J)Y
=-r(JY). \eea

\subsection{The canonical connection and some of its properties}
An important property of nearly K\"ahler geometry is the existence of a
canonical Hermitian connection \cite{FI} (see \cite{CS} for pseudo-Riemannian metrics). This is the unique connection $\bar \n$ with skew-symmetric 
torsion, which parallelizes the metric $g$ and the almost complex structure $J.$ Explicitely it is given
by 
\be \bar \n_X Y = \n_X Y -\frac{1}{2} J(\n_X J)Y, \mbox{ for } X, Y \in
\Gamma(TM) \label{can_con_def} \ee
and its torsion equals $\bar{T}(X,Y)= -J(\n_X J)Y.$
\bp \label{parallel_torsion}
Denote by $\bar \n$ the canonical connection of a nearly pseudo-K\"ahler manifold $(M,J,g).$ 
Then one has $\bar\n(T)=\bar\n (\n J)=0.$ 
\ep
\begin{proof} The proof given in \cite{BM} essentially uses the
explicit form \eqref{can_con_def} of the connection $\bar\n$ and the identity
\eqref{id_equ_D2_J}. Therefore the proposition generalizes to the
pseudo-Riemannian case.
\end{proof}
 From Proposition \ref{parallel_torsion} and the relation \eqref{can_con_def} of $\n$ and $\bar \n$ one obtains the
 following identities for the curvature tensor $\bar R$ of $\bar \n$ and the curvature tensor $R$ of the Levi-Civita connection $\n$
\bea 
\bar{R}(W,X,Y,Z)&=&  {R}(W,X,Y,Z)-\frac{1}{2}g((\n_WJ)X,(\n_YJ)Z)\nonumber \\
&+& \frac{1}{4}\left[g((\n_WJ)Y,(\n_{X}J)Z)-
  g((\n_WJ)Z,(\n_{X}J)Y)\right]  \label{curv_rel_n_nbar}\\ 
&=& \frac{1}{4} [ 3R(W,X,Y,Z) + R(W,X,JY,JZ) \nonumber \\ &+&\sigma_{XYZ}
  R(W,X,JY,JZ)],\nonumber \\
\bar{R}(W,JW,Y,JZ)&=&\frac{1}{4}[ 5 R(W,JW,Y,JZ)\nonumber \\&-&
 R(W,Y,W,Z)  -R(W,JY,W,JZ)]. \label{curv_rel_n_nbar_prime}
\eea
With the help of the  equation \eqref{curv_rel_n_nbar} it follows
\bea \quad \quad
\bar R(W,X,Y,Z)= \bar R(Y,Z,W,X)= -\bar R(X,W,Y,Z) = -\bar R(W,X,Z,Y).\label{curiv_id_i}  \eea
Using $\bar \n J=0$ and $\bar \n g=0$ we obtain
\bea
\bar R(W,X,Y,Z)&=& \bar R(W,X,JY,JZ) \label{cx_curv_id} \\ \nonumber &=& \bar R(JW,JX,Y,Z)= \bar R(JW,JX,JY,JZ).
\eea
The general form of the first Bianchi identity (cf. chapter III of
\cite{KN}) for a connection with torsion
yields in the case of parallel torsion:
\bea
\underset{XYZ}{\sigma}  \, \bar R(W,X,Y,Z) &=& -\underset{XYZ}{\sigma} \, g( (\n_WJ)X,
(\n_YJ)Z). \label{curiv_id_iii}
\eea
In a similar way we get from the second Bianchi identity
(cf. chapter III of \cite{KN}) for a connection
with parallel torsion or from the second Bianchi identity for $\n$ 
\bea 
-\underset{VWX}\sigma \bar \n_V (\bar R)(W,X,Y,Z) =  \underset{VWX}\sigma \bar R((\n_VJ)JW,X,Y,Z). \label{curiv_id_bianci_ii}
\eea
From deriving  equation \eqref{curv_rel_n_nbar} and the second Bianchi identity of $\n$
one gets after a direct computation
\bea
\underset{VWX}\sigma \n_V (\bar R)(W,X,Y,Z) 
= \frac{1}{2} g((\n_Y)Z,\underset{VWX}\sigma (\n_XJ)(\n_VJ)JW) \label{curiv_id_iv},
\eea
which implies
\bea
\underset{VWX}\sigma \n_V (\bar R)(W,X,Y,JY) =0. \label{curiv_id_iv_prime}
\eea
\bp
The tensor $r$ on a nearly pseudo-K\"ahler manifold $(M,J,g)$ is parallel with respect to the canonical connection $\bar \n.$
\ep
\begin{proof}
Deriving $g(rX,X) = \sum_{i=1}^{2n} \epsilon_i g((\n_{e_i}J)X,(\n_{e_i}J)X)$ one
obtains
\bean
g((\n_Ur)X,X)&=& 2 \sum_{i=1}^{2n} \epsilon_i g((\n^2_{U,e_i}J)X,(\n_{e_i}J)X) \\
&\overset{\eqref{id_equ_D2_J}}=&  -\sum_{i=1}^{2n}  \epsilon_i
[g((\n_UJ)e_i, (\n_XJ) J (\n_{e_i}J)X ) \\
&+& g((\n_UJ) (\n_{e_i}J)X, (\n_{e_i}J)JX) \\
&+&  g((\n_UJ)X, (\n_{(\n_{e_i}J)X} J) J e_i)].
\eean
We observe that in the first two terms exchanging $e_i$ by $J e_i$ gives a
minus sign. Hence taking an adapted orthogonal frame $\{e_i\}_{i=1}^{2n}$ yields:
\bean
g((\n_Ur)X,X)&=&   \sum_{i=1}^{2n}  \epsilon_i g((\n_UJ)X, (\n_{J e_i}
J)(\n_{e_i}J)X ) =   -\sum_{i=1}^{2n}  \epsilon_i g((\n_UJ)X, (\n_{ e_i} J)(\n_{e_i}J)JX )\\
&=&\sum_{i=1}^{2n}  \epsilon_i g((\n_{ e_i}
J)(\n_UJ)X, (\n_{e_i}J)JX )=g(r(\n_UJ)X,JX).
\eean
Polarizing this expression shows using that $r$ is $g$-symmetric the
following identity
$$g\left((\n_Ur)X,Y\right) = \frac{1}{2} g(r(\n_UJ)X,JY)+  \frac{1}{2}
g(JX,r(\n_UJ)Y).$$
As the difference of $\bar \n$ and $ \n$ is  $-\frac{1}{2}J\n J$ the
last equation is exactly $\bar \n r =0.$
\end{proof}

\bt \label{THM_curv} Let $(M,J,g)$ be a nearly pseudo-K\"ahler manifold and let $W,X$ be vector fields on $M$ then it holds
\be \sum_{i,j=1}^{2n} \epsilon_i \epsilon_j
g(re_i,e_j)
\left[R(W,e_i,X,e_j) -5 R(W,e_i,JX,Je_j) \right] =0.\ee
\et
\begin{proof} Let $\{ e_i\}_{i=1}^{2n}$ be an adapted 
  orthogonal frame field. One observes
\bean
 \sum_{i=1}^{2n} \epsilon_i 
\bar R(W,X,e_i,(\n_VJ)e_i) \overset{\eqref{cx_curv_id}}= 
~\frac{1}{2}  \sum_{i=1}^{2n} &\epsilon_i
\left[\bar R(W,X,e_i,(\n_VJ)e_i)\right. & \\ 
&- \left. \bar R(W,X,Je_i,(\n_VJ)Je_i)\right]&=0
\eean
and one gets after derivation of the left hand-side
\be \sum_{i=1}^{2n}  \epsilon_i\left[\n_U\left(\bar R)(W,X,e_i,(\n_VJ)e_i \right)  
+ \bar R\left(W,X,e_i,\n_{U,V}^2(J)e_i \right)\right]=0. \label{zw_eq} \ee
Taking the trace in $U,V$ and applying the identity \eqref{id_equ_spur_D2_J} yields on the second term  
\bean &~&\sum_{i,k=1}^{2n}  \epsilon_i\epsilon_k\left[\bar
  R\left(W,X,e_i,\n_{e_k,e_k}^2(J)e_i  \right)\right] = 
\sum_{i=1}^{2n}  \epsilon_i \bar
  R\left(W,X,e_i,-r(Je_i)  \right)
\\  &=&-\sum_{i,j=1}^{2n}  \epsilon_i \epsilon_j
g(re_i,e_j)\bar R\left(W,X,e_i,Je_j   \right).
\eean
From \eqref{curv_rel_n_nbar_prime} we obtain
\bean 4\bar{R}\left(W,JW,e_i,Je_j\right)= 
5R(W,JW,e_i,Je_j)-R(W,e_i,W,e_j)-R(W,Je_i,W,Je_j).
\eean
The first Bianchi identity using an adapted frame implies for $i,j \in \{ 1,\ldots,n\}$:
\bean R(W,JW,e_i,Je_j)&=&  -R(JW,e_i,W,Je_j)  -R(e_i,W,JW,Je_j)  \\
&=& -R(W,Je_j,JW,e_i)  +R(W,e_i,JW,Je_j)\\
&=& R(W,e_{j+n},JW,Je_{i+n})  +R(W,e_i,JW,Je_{j}).
\eean
In an adapted frame it is $g(re_{i+n},e_{j+n}) = g(re_i,e_j)$ with $i,j \in \{
1,\ldots,n\}.$ Therefore taking the trace in an adapted frame and polarizing yields the claimed
identity if we show the vanishing of the trace on the first  term of \eqref{zw_eq}.
\bean
 \sum_{i,j=1}^{2n}  \epsilon_i \epsilon_j\n_{e_j}\left(\bar R)(W,JW,e_i,(\n_{e_j}J)e_i \right) 
&=& \sum_{i,j,k=1}^{2n}  \epsilon_i \epsilon_j\epsilon_k \, g((\n_{e_j}J)e_i,e_k
 )\, \n_{e_j}\left( \bar R \right)(W,JW,e_i,e_k )
 \\
&=& \frac{1}{3} \sum_{i,j,k=1}^{2n} \underset{ijk}\sigma \left[
\epsilon_i \epsilon_j \epsilon_k g((\n_{e_j}J)e_i,e_k ) \n_{e_j}(\bar R)(W,JW,e_i,e_k )\right].
\eean
Since $\epsilon_i \epsilon_j \epsilon_k g((\n_{e_j}J)e_i,e_k )$ is
  constant under cyclic permutation of $i,j,k,$ the last expression vanishes as a
  consequence of the curvature identities \eqref{curiv_id_i} and the identity \eqref{curiv_id_iv_prime}. Polarization in $W$
finishes the proof.
\end{proof}
\section{First structure results}
\subsection{Small dimensions} \label{small_dim_sect}
\noindent For a nearly pseudo-K\"ahler manifold $\n \o$ is a differential form of type
$(3,0)+(0,3).$ In consequence real two- or four-dimensional nearly pseudo-K\"ahler manifolds are automatically pseudo-K\"ahler.
Six dimensional nearly pseudo-K\"ahler manifolds are either   pseudo-K\"ahler manifolds or strict
nearly pseudo-K\"ahler manifolds. In the strict case a nearly pseudo-K\"ahler manifold $(M^6,J,g)$ is of {\cmssl constant type}, i.e. it holds
\be 
g((\n_XJ)Y,(\n_XJ)Y)= \alpha \left(g(X,X)g(Y,Y)- g(X,Y)^2 -g(JX,Y)^2\right).
\ee
The sign of the type constant $\alpha$ depends on the signature $(p,q)$ with $p+q=6$.
In fact it is $\mbox{sign}(\alpha)= \mbox{sign}(p-q),$ see section 7 of
\cite{Ka}.
\subsection{Linear algebra of three-forms}
In the following section we consider a (finite dimensional) pseudo-hermitian vector 
space $(V,J,\langle \cdot, \cdot \rangle).$ Let $\eta\in \Lambda^3V^*$ 
 be a three-form. We define the {\cmssl support} of $\eta$ by 
\bea \Sigma_{\eta} = \mbox{span} \{ \, X \lrcorner Y \lrcorner \eta \, |\, X,Y \in V \, \} \subset
V,\eea
where we identified $V$ and $V^*$ by means of $\langle \cdot, \cdot \rangle.$
The name support is motivated by the observation, that for a given 
$\eta \in \Lambda^3V^*$ it already holds $\eta \in \Lambda^3 \Sigma_\eta^*,$
 compare Lemma 7 of \cite{CS}.

In the present paper we are essentially interested in the three-form 
$g_p((\n_XJ)_pY,Z)$ for $X,Y,Z \in T_pM$ on a nearly pseudo-K\"ahler manifold
$(M,J,g).$ This three-form is a real form of type $(3,0)+(0,3).$ The type
condition implies that $\Sigma_{\eta}$ is a $J$-invariant subspace. In
particular it follows that the complex dimension of the support of a non-zero
 such form is at least three.\\
In \cite{CS} the classification of Levi-Civita flat nearly K\"ahler manifolds
was related to the existence of real three-forms of type $(3,0)+(0,3)$
 with isotropic  support, i.e. such that  $\Sigma_\eta$ is an
isotropic subspace. \\
We define the kernel of a three-form $\eta \in \Lambda^3V^*$ by $ \mathcal K= \mathcal K_\eta = \mbox{ker}(X \mapsto X  \lrcorner \eta).$

\bl
 One has $\mathcal{K}= \Sigma_\eta^\bot$ and  $\Sigma_\eta =
 \mathcal{K}^\bot.$ 
\el
\begin{proof}
Suppose, that $X$ is in $\mathcal K_\eta,$ i.e. $X \lrcorner \eta=0.$
By definition $\Sigma_\eta$ is spanned by vectors $U$ satisfying
$\langle U,\cdot \rangle =\eta(Y,Z,\cdot)$ for $Y,Z \in V.$ This implies
$\langle U,X \rangle=\eta(Y,Z,X)=0,$ since $\eta$ is a
three-form. If $X$ is perpendicular to $\Sigma_\eta$ the claim follows from
the last equation and  $\eta \in\Lambda^3\Sigma_\eta^*.$ This means  $X$ is in $\mathcal K_\eta$ if and only if $X$ is
perpendicular to $\Sigma_\eta.$ It follows $\Sigma_\eta = \mathcal{K}^\bot.$ 
\end{proof}

\bl \label{four_dim_form_lemma}
Let $(V,J,\langle \cdot,\cdot \rangle )$ be a pseudo-hermitian vector space with $\mbox{dim}_\bR(V)=8$ 
then a real three-form $\eta$ of  type $(3,0)+(0,3)$ and of non-vanishing length 
has a (complex) one dimensional kernel $\mathcal K_\eta,$ which admits an orthogonal complement  $(\mathcal
K_\eta)^\bot.$ Moreover one has $\Sigma_\eta=(\mathcal K_\eta)^\bot.$
\el
\begin{proof} Let us identify $V$ and $V^{1,0}.$ Denote by $\{\alpha_i \}_{i=1}^4$ a unitary basis of $(V^{1,0})^*$
and define  a $(4,0)$-form $v$ by $v= \alpha_1\wedge \alpha_2\wedge\alpha_3 \wedge \alpha_4.$ 
The map 
 \bean \Phi \,:\, V^{1,0}\ra \Lambda^3(V^{1,0})^*, \quad
 \zeta \mapsto \zeta \lrcorner  v
\eean
 yields an isometry. Therefore the $(3,0)$-form $\rho= \eta + i J^*\eta$ 
is given by $\rho = Z \lrcorner  v$ for some  $Z\in V^{1,0}$ 
and consequently it follows $Z \lrcorner \rho=0.$ 
As $\Phi$ is an isometry and $\rho$ has non-zero length, we conclude 
that $Z$ is not isotropic. Denote by $\mathcal L \subset \mathcal K_\rho$ the complex line 
spanned by $Z$  and by $\mathcal{L}^\bot\supset \mathcal{K}^\bot_\rho$ its orthogonal  complement. It
remains to prove $\mathcal L= \mathcal K:$ On the one hand we have  $\Sigma_\rho = \mathcal{K}^\bot_\rho\subset\mathcal{L}^\bot$ and on the
other hand from $\rho\ne 0$ we get
$\mbox{dim}_\bC \Sigma_\rho \ge 3 = \mbox{dim}_\bC \mathcal{L}^\bot.$ 
This shows  $\Sigma_\rho = \mathcal{L}^\bot$ and $\mathcal K_\rho = \mathcal L.$
\end{proof}
\begin{remark}
As the reader observes, if $\eta$ has length zero, one can replace the
orthogonal complement by the null-space and obtain an analogous statement as in the last proposition.
\end{remark}
\bl \label{threeform_dim10}
Let $(V,J,\langle \cdot,\cdot \rangle )$ be a pseudo-hermitian vector space with $\mbox{dim}_\bR(V)=10$ 
then a  real three-form $\eta$ of  type $(3,0)+(0,3)$ and of non-vanishing
length is of the following possible types:
\begin{itemize}
\item[(i)] There exists an orthonormal real basis $\{f_i\}_{i=1}^{10}=\{ e_1,Je_1,\ldots , e_5,
  Je_5\}$  and real numbers $\alpha,\beta$ such that
\bea
\eta(e_1,e_2,e_3)= \alpha \ne 0; \quad \eta(e_4,e_5,e_1)= \beta \label{norm_form_10_1} 
\eea
and $\eta(f_i,f_j,f_k)=0$ for the cases which are not obtained
from \eqref{norm_form_10_1} by skew-symmetry and type relations.
\item[(ii)] There exists an orthonormal real basis $\{f_i\}_{i=1}^{10}=\{ e_1,Je_1,\ldots , e_5,
  Je_5\}$  and real numbers $\alpha,\beta$ such that
\bea
\eta(e_1,e_2,e_3)= \alpha \ne 0; \quad \eta(e_4,e_5,e_1+e_3)= \beta
\mbox{ with } \langle e_1,e_1 \rangle =-\langle e_3,e_3 \rangle \label{norm_form_10_2} 
\eea
and $\eta(f_i,f_j,f_k)=0$ for the cases which are not obtained
from \eqref{norm_form_10_2} by skew-symmetry and type relations.
\end{itemize}
\el
\begin{proof}
Denote by $\{\alpha_i \}_{i=1}^5$ a unitary basis of $(V^{1,0})^*$
and define  a $(5,0)$-form $v$ by $v= \alpha_1\wedge \alpha_2\wedge\alpha_3 \wedge \alpha_4\wedge \alpha_5.$ 
The map 
 \bean \Phi \,:\, \Lambda^2V^{1,0}\ra \Lambda^3(V^{1,0})^*, \quad
 \varphi \mapsto \varphi \lrcorner  v
\eean
 yields an isometry. Therefore the $(3,0)$-form $\rho= \eta + i J^*\eta$ 
is given by $\rho = \varphi \lrcorner  v$ for some  $\varphi \in \Lambda^2 V^{1,0}.$ As $\Phi$ is an isometry and $\eta$ has non-zero length, we conclude 
that $\varphi$ is not isotropic. Define $Z \in (V^{1,0})^*$ by $Z=\varphi
\lrcorner \rho=\varphi \lrcorner (\varphi \lrcorner \rho).$ From 
$\langle Z,Z\rangle =\langle\varphi,\varphi\rangle^2 \langle \rho,\rho\rangle$ we obtain that $Z$ is not isotropic.\\
Choosing a unitary basis
$\{\frac{1}{\sqrt{2}}(e_1-iJe_1),\frac{1}{\sqrt{2}}(e_2-iJe_2)\}$ of the plane 
  $\varphi$ and $Z= \alpha' \frac{1}{\sqrt{2}}(e_3 - iJe_3)^*$ for a unit vector
 $e_3\in V$ and $\alpha' \in \bR-\{0\}$ we consider
$\mathcal B^{1,0}:= \{ \frac{1}{\sqrt{2}}(e_1-iJe_1),\frac{1}{\sqrt{2}}( e_2-iJe_2),\frac{1}{\sqrt{2}}(e_3 -iJe_3) \}.$\\
{\bf Claim:} For $\zeta,\chi  \in \mbox{span}_\bC \, \mathcal{B}^{1,0}$ it follows
 $\rho(\zeta,\chi,\cdot)  \in \mbox{span}_\bC \, (\mathcal{B}^{1,0})^*.$ \\
Let us  define the map
$
 \tilde \Phi \,:\, \Lambda^2(\mbox{span}_\bC \, \mathcal{B}^{1,0}) \ra
 (V^{1,0})^*,$ by linear extension of 
$
\zeta \wedge \chi \mapsto \rho(\zeta,\chi,\cdot).$ \\
We observe that $\Lambda^2(\mbox{span}_\bC\, \mathcal{B}^{1,0})$ is a vector space 
of complex dimension three, which implies that $\dim_\bC(\mbox{im}\,\tilde\Phi)
\le 3.$ As $\rho$ is a three-form one easily sees that the duals of $ e_j
-iJe_j$ for $j=1,2,3$ are contained in $\mbox{im}\,\tilde \Phi.$ By the bound
on the dimension of $ \mbox{im}\,\tilde \Phi$ the components orthogonal to $ (\mbox{span} \,\mathcal{B}^{1,0})^*$ vanish. This proves the claim.\\
Let $W$ be the orthogonal complement of $\mbox{span}_\bC \mathcal B^{1,0}.$ 
Choose an orthogonal basis $\{ e_4 -iJe_4, e_5 -iJe_5\}$ of $W.$ Using that
 $\rho$ is skew-symmetric we conclude that  $\tilde Z=\rho( e_4 -iJe_4, e_5
-iJe_5,\cdot)$ is perpendicular to the dual $W^*$ of $W$ and hence an element
$\tilde Z$ of $(\mbox{span}_\bC \mathcal B^{1,0})^*.$ 
If $\langle \tilde Z, \tilde Z \rangle \ne 0$ we can adapt the basis of $\mathcal B^{1,0}$
 such that $\tilde Z= \beta' \frac{1}{\sqrt{2}}(e_1-iJe_1)^*.$ If 
$\langle \tilde Z,\tilde Z \rangle =0$ we can achieve  $\tilde Z= \beta' \frac{1}{\sqrt{2}}
\left[(e_1 -iJe_1) + (e_3 -iJe_3)\right]$ with  $\langle e_1,e_1 \rangle=-
\langle e_3,e_3 \rangle.$ Passing to the real basis yields some new constants
$\alpha, \beta$ and the claim of the Lemma.
\end{proof}

\subsection{K\"ahler factors and the structure in dimension 8}
The aim of this subsection is to split off the pseudo-K\"ahler factor of a
nearly pseudo-K\"ahler manifold. This will be done by means
of the kernel of $\n J$ and allows to reduce the (real) dimension from eight
to six. \\
For $p\in M$ we set  
$$ \mathcal K_p = \mbox{ker}(X \in T_pM \mapsto \n_X J).$$
\bt \label{Kaehler_factor_prop} Let $(M,J,g)$ be a nearly pseudo-K\"ahler
manifold.  Suppose, that the distribution $\mathcal K$ has constant dimension
and admits an orthogonal complement,
\begin{itemize}
\item[(i)] then $M$ is locally a pseudo-Riemannian product $M= K \times M_1$
  of a pseudo-K\"ahler manifold $K$ and a strict nearly pseudo-K\"ahler manifold $M_1.$
\item[(ii)] if $M$ is complete and simply connected then it is a pseudo-Riemannian product $M= K \times M_1$ of a
pseudo-K\"ahler manifold $K$ and a strict nearly pseudo-K\"ahler manifold $M_1.$
\end{itemize}
\et
\begin{proof}
The distribution $\mathcal K$ is parallel for the canonical connection $\bar
\n,$ since $\n J$ is $\bar \n$-parallel. By the formula \eqref{can_con_def}
and the nearly K\"ahler condition it follows $\bar \n_XK = \n_XK$ for sections $K$ in $\mathcal K$ and 
$X$ in $TM.$ This implies that  $\mathcal K$ is parallel for the Levi-Civita
connection and in consequence its orthogonal complement $(\mathcal K)^\bot$ 
is Levi-Civita parallel. The proof of (i) finishes by the local version of the theorem of de Rham and
the proof of (ii) by the global version.
\end{proof}

\begin{remark}
There exist nearly pseudo-K\"ahler manifolds $(M,J,g)$ without pseudo-K\"ahler
de Rham factor, such that $\mathcal K_\eta\ne \{ 0\}$ admits no orthogonal
complement. In fact  there are Levi-Civita flat nearly pseudo-K\"ahler manifolds constructed in
Theorem 1 and 3 of \cite{CS} such that the three-form
$\eta_p(X,Y,Z)=g_p(J(\n_XJ)Y,Z),$ for  $p\in M,$ has a support $\Sigma_\eta \subset T_pM$ which is a maximally isotropic subspace~(Here we identified $T_pM$ and $T^*_pM$ via the
metric $g.$). Obviously, $J(\n_XJ)Y$ and  $J(\n_UJ)V$ are elements of the
support of $\eta$ for arbitrary $X,Y,U,V\in T_pM.$ It then follows
$0=g(J(\n_XJ)Y,J(\n_UJ)V)=g(J(\n_{J(\n_XJ)Y}J)U,V)$ for all $V\in T_pM.$ 
Hence it is  $\Sigma_\eta \subset \mathcal K_\eta.$ Moreover
for general reasons we have shown before $\Sigma_\eta = \,\mathcal{K}_\eta^\bot$ which shows $\mathcal K_\eta\cap  \mathcal K_\eta^\bot \ne \{ 0\}$ for the
above examples. From these examples we learn, that the Theorem \ref{Kaehler_factor_prop}
does not hold true, if there is no orthogonal complement. 
\end{remark}

\bd \label{nice_def}
A nearly pseudo-K\"ahler manifold  $(M,J,g)$ is called {\cmssl nice} if 
 the three-form $g((\n_\cdot J)\cdot,\cdot)$ has non-zero length in each point
$p\in M.$
\ed
\bt \label{class_dim8}
Let $(M^8,J,g)$ be a complete simply connected eight-dimensional nice nearly pseudo-K\"ahler
manifold. Then $M = M_1 \times M_2$ where $M_1$ is a
two-dimensional K\"ahler manifold and $M_2$ is a six-dimensional strict nearly
pseudo-K\"ahler manifold.
\et
\begin{proof}
Since $(M,J,g)$ is a nice nearly pseudo-K\"ahler manifold we can use Lemma \ref{four_dim_form_lemma}
 to obtain an orthogonal splitting in the two-dimensional distribution 
$\mathcal K$ and its  orthogonal complement, which coincides with
 $\Sigma_\eta.$ Therefore we are in the situation of Theorem
 \ref{Kaehler_factor_prop} (ii).
\end{proof}
\subsection{Einstein condition versus reducible holonomy}
\bt \label{Einst_theo}
Let $(M,J,g)$ be a nearly pseudo-K\"ahler manifold.
\begin{itemize}
\item[(i)]
Suppose that $r$ has more than one eigenvalue, then the canonical Hermitian 
connection has reduced holonomy.
\item[(ii)]
If the tensor field $r$ has exactly one eigenvalue then $M$ is a
pseudo-Riemannian Einstein manifold.
\end{itemize}
\et
\begin{proof}
(i) Let $\mu_i$ for $i=1, \ldots, l$ be the eigenvalues of $r.$ Then the
  decomposition in the according eigenbundles $\mbox{Eig}(\mu_i)$ is $\bar \n$-parallel and hence
its holonomy is reducible. \\
(ii) From the identity of Theorem \ref{THM_curv} and $r = \mu
\id_{TM}$ we obtain 
$$ 0=\sum_{i=1}^{2n} \epsilon_i 
\left(R(W,e_i,X,e_i) -5 R(W,e_i,JX,Je_i) \right) =g((Ric-5Ric^*)W,X),$$ 
where we used the Bianchi identity and an adapted frame to obtain the last equality. This shows comparing with $r=Ric-Ric^*$ that it holds $Ric= \frac{5}{4} \mu.$ 
\end{proof}

Let us recall, that in the pseudo-Riemannian setting the decomposition into the eigenbundles 
is {\bf not} automatically ensured to be an orthogonal direct decomposition. Therefore we introduce the following notion:
\bd \label{decomp_def}
A  nearly pseudo-K\"ahler manifold $(M,J,g)$ is 
called {\cmssl decomposable} if the above decomposition into the eigenbundles of the tensor $r$ is orthogonal.
\ed

\bt \label{Einst_dim_6}
A strict nearly pseudo-K\"ahler six-manifold $(M^6,J,g)$ of constant type
$\alpha$ is a pseudo-Riemannian Einstein manifold with Einstein constant $5\alpha$.
\et
\begin{proof} In an adapted basis we obtain from the symmetries of $\n J$
\bean
g(rX,X)&=&2\sum_{i=1}^{3} \epsilon_i\, g((\n_{X}J)e_i, (\n_{X} J)e_i )
=-2\sum_{i=1}^{3} \epsilon_i\, g((\n_{X}J)^2e_i, e_i).\eean
This  is exactly minus the trace of the operator $(\n_XJ)^2$ which has a simple form
in a cyclic frame. It follows after polarizing $g(rX,Y)=4\alpha g(X,Y).$ From 
Theorem \ref{Einst_theo} we compute the Einstein constant $5 \alpha$ where
$\alpha$ is the type constant of the strict nearly pseudo-K\"ahler manifold $M^6.$
\end{proof}

\bp \label{berechne_r}
Let $(M^{10},J,g)$ be a nice nearly pseudo-K\"ahler ten-manifold.
\begin{itemize}
\item[(i)]
Then the tensor $r$ in a frame of the first type in Lemma \ref{threeform_dim10} is given by
\bean
r e_1&=& 4( \alpha^2 + \beta^2)e_1, \\
r e_2&=& 4 \alpha^2 e_2, \quad 
r e_3= 4 \alpha^2 e_3, \\
r e_4&=& 4 \beta^2 e_4, \quad
r e_5= 4 \beta^2 e_5, \\
r(Je_i)&=& Jr(e_i), \quad i=1,\ldots, 5,
\eean
where $\alpha,\beta$ are constants.
\item[(ii)]
For a frame of the second type in Lemma \ref{threeform_dim10} the tensor $r$ is given by

\[ r\left[ \begin {array}{c} e_1 \\ e_2 \\ e_3
\end{array} \right]
=4\left[ \begin {array}{ccc}
    {\alpha}^{2}+\,{\beta}^{2}\epsilon_{{4}}\epsilon_{{5}}&0&\,{\beta}^{2}\epsilon_{{4}}\epsilon_{{5}}\\
\noalign{\medskip}0&\,{\alpha}^{2}&0\\
\noalign{\medskip}\,{\beta}^{2}\epsilon_{{4}}\epsilon_{{5}}&0&\,{\alpha}^{2}+\,{\beta}^{2}\epsilon_{{4}}\epsilon_{{5}}
\end{array} \right]\left[ \begin {array}{c} e_1 \\ e_2 \\ e_3
\end{array} \right] \]
\bean
r e_4&=& 0, \\
r e_5&=& 4\,{\beta}^{2}(2\epsilon_1\epsilon_4 -1)e_5,\\
r(Je_i)&=& Jr(e_i), \quad i=1,\ldots, 5.
\eean
The eigenvalues are   $\{0;4{\alpha}^{2};4\,{\beta}^{2}(2\epsilon_1\epsilon_4 -1);4(\alpha^2+2\beta^2 \epsilon_4\epsilon_5)\}, $  where the eigenbundles are given as
\bean
{\rm Ker}(r)&=&{\rm span}\{e_4,J e_4\}, \\
{\rm Eig}(r,4\alpha^2)&=&{\rm span}\{ -e_1+e_3,e_2,-Je_1+Je_3,Je_2\}, \\
{\rm Eig}(r,4\,{\beta}^{2}(2\epsilon_1\epsilon_4 -1))&=&{\rm span}\{e_5,Je_5\},\\
{\rm Eig}(r,4(\alpha^2 +2\beta^2 \epsilon_4\epsilon_5)))&=&{\rm span}\{e_1+e_3,Je_1+Je_3\},
\eean
where $\alpha,\beta$ are constants. For $\beta^2 \ne 0$ the second case is not decomposable. 
\item[(iii)] Suppose $\beta =0$ in the case (i) and (ii). Then it follows
\bean
{\rm Eig}(r,4\alpha^2)&=&{\rm span}\{ e_1,e_2,e_3,Je_1,Je_2,Je_3 \},\\
{\rm Ker}(r)&=&{\rm span}\{ e_4,e_5,Je_4,Je_5 \}.
\eean
\end{itemize}
\ep
\begin{proof}In an adapted basis we obtain from the symmetries of $\n J$
\bean
g(rX,Y)&=&2\sum_{i=1}^{5} \epsilon_i\, g((\n_{X}J)e_i, (\n_{Y} J)e_i )
=-2\sum_{i=1}^{5} \epsilon_i\, g((\n_YJ)(\n_XJ)e_i, e_i).\eean
This is exactly minus the trace of the operator $(\n_YJ)(\n_XJ).$ Using the form of Lemma \ref{threeform_dim10} one can calculate $r$ by hand or using computer
algebra systems to obtain the claimed results.
\end{proof}

\bt \label{split_dim_ten_thm}
Let $(M^{10},J,g)$ be a complete simply connected nice decomposable nearly
pseudo-K\"ahler manifold of dimension ten. Then $M^{10}$ is of one of the following types
\begin{itemize}
\item[(i)]
the tensor $r$ has a kernel and $M^{10}= K \times M^6$ is a product of a four-dimensional
pseudo-K\"ahler manifold $K$ and a strict nearly pseudo-K\"ahler six-manifold $M^6.$
\item[(ii)] the tensor $r$ has trivial kernel and $r$ has eigenvalues
  $4(\alpha^2+\beta^2)$ with multiplicity $2,$ $4\alpha^2, 4\beta^2$ with
  multiplicity $4$ for some $\alpha, \beta \ne 0,$
\end{itemize}
A nice nearly pseudo-K\"ahler manifold $(M^{10},J,g)$ is decomposable if the
dimension of the kernel of $r$ is not equal to two.
\et
\begin{proof}
Since we suppose, that $(M^{10},J,g)$ is a nice and decomposable nearly
pseudo-K\"ahler manifold, Proposition \ref{berechne_r} implies that one has
the two different cases:\\
(i) the distribution $\mathcal K,$ which is the tangent
  space of the K\"ahler factor has dimension four  and admits an orthogonal
  complement of dimension six. This is part (iii) of  Proposition
  \ref{berechne_r}. Part (i) of the Theorem now follows from Theorem \ref{Kaehler_factor_prop}.\\
(ii) the tensor $r$ has trivial kernel and we are in the situation of
  Proposition  \ref{berechne_r} part (i) with $\alpha,\beta \ne 0$ and part (ii) follows.
\end{proof}

\begin{remark}
Nearly pseudo-K\"ahler manifolds falling in the second case of the last
theorem will be shown to be related to twistor spaces in section \ref{tw_proof}.
\end{remark}

\section{Pseudo-Riemannian submersions}
Let us consider  the setting of a pseudo-Riemannian submersion 
$\pi :(M,g) \ra (N,h).$ The tangent bundle $TM$ of $M$ splits orthogonally 
into the direct sum \bea TM= \mathcal{H} \oplus \mathcal{V}. \label{split_TM} \eea
Denote by $\iota_{\mathcal{H}} ,\iota_{\mathcal{V}}$ the canonical inclusions 
and by $\pi_{\mathcal{H}} ,\pi_{\mathcal{V}}$ the canonical projections. We recall
the definition \cite{Besse,ON} of the fundamental tensorial invariants $A$ and
$T$ of the submersion $\pi$ 
\bean
T_{\zeta}  = \pi_{\mathcal{H}}\circ \nabla_{\pi_{\mathcal{V}}(\zeta)}\circ \pi_{\mathcal{V}} + \pi_{\mathcal{V}}\circ \nabla_{\pi_{\mathcal{V}}(\zeta)}\circ\pi_{\mathcal{H}},
  \\   
A_\zeta  = \pi_{\mathcal{H}}\circ \nabla_{\pi_{\mathcal{H}}(\zeta)}\circ \pi_{\mathcal{V}} + \pi_{\mathcal{V}}\circ \nabla_{\pi_{\mathcal{H}}(\zeta)}\circ\pi_{\mathcal{H}},
\eean
where $\zeta$ is a vector field on $M.$ \\
The components of the Levi-Civita connection $\n$ are given in the next
proposition (compare \cite{ON}, \cite{Besse} 9.24 and 9.25).
\bp
Let $\pi :(M,g) \ra (N,h)$ be a pseudo-Riemannian submersion, denote by
$\n$  the Levi-Civita connection of  $g$ and define ${\n}^{\mathcal{V}}:=
\pi_{\mathcal{V}} \circ \n \circ \iota_{\mathcal{V}}.$ For vector fields $X,Y$ in
$\mathcal{H}$ and $U,V$ in $\mathcal{V}$ we have the following
identities 
\bea
\n_UV&=& {\n}^{\mathcal{V}}_UV +T_UV, \label{sub_zush_equ_I}\\
\n_UX&=&T_UX+ \pi_{\mathcal{H}}(\n_UX),\\
\n_XU&=& \pi_{\mathcal{V}}(\n_XU)+ A_XU,\\
\n_XY&=& A_XY + \pi_{\mathcal{H}}(\n_XY),\\
\pi_{\mathcal{V}} [X,Y] &=& 2A_XY, \label{A_and_Lie}
\\
g(A_XY,U) &=&-g(A_XU,Y),  \mbox{ or more generally } A \mbox{ is alternating.}  \label{A_altern}
\eea
\ep
\noindent
The {\cmssl canonical variation} of the metric $g$ for $t\in \bR -\{0\}$ is given by
$$ {g}_t:= 
\begin{cases} 
 g(X,Y), \mbox{ for } X,Y \in \mathcal{H}, \\
 t g(V,W), \mbox{ for } V,W \in \mathcal{V},\\
g(V,X)=0,   \mbox{ for } V \in \mathcal{V},
X \in \mathcal{H} .
\end{cases}$$
\bl \label{Besse_Lemma}
Denote by $X,Y$ vector fields in  $\mathcal H$ and by $U,V$  vector fields in $\mathcal V.$
\begin{enumerate}
\item
 Let $A^t$ and $T^t$ be the tensorial invariants for $g_t$ and $A$ and $T$
 those for $g=g_1.$  Then it holds 
\bea \label{Besse_Lemma_I_equ}&~&A^t_XY=A_XY,\quad A^t_XU=tA_XU \mbox{  and  } \\
&~& T^t_UV=tT_UV,\quad T^t_UX=T_UX; \label{Besse_Lemma_Iprime_equ}\eea
\item 
 $\n_U^{t\,\mathcal{V}}V={\n}^{\mathcal{V}}_UV;$
\item
\bea \label{Besse_Lemma_III_equ}\pi_{\mathcal{H}} (\n^t_X Y) =
  \pi_{\mathcal{H}} ({\n}_X Y)
 \mbox{ and } \pi_\mathcal{V}(\n^t_XV)= \pi_\mathcal{V}({\n}_XV); \eea
\item \bea &~&\,\pi_{\mathcal{H}}( {\n^t}_V X  -\n_V X) =  (t-1)A_X
  V; \label{Diff_n_nz_Hpart} \\
&~&\n^t_UV=\n_UV+(t-1)T_UV.
\eea
\end{enumerate}
\el
\begin{proof}
The first part can be found in Lemma 9.69 of \cite{Besse}.
On the right hand-side of the Koszul formulas one only needs 
the metric $g_t$ on $\mathcal{V}$ to determine ${{\n}^t}^{\mathcal{V}}.$ This shows ${{\n}^t}^{\mathcal{V}}=\n^{\mathcal{V}}.$ An analogous argument using the Koszul formulas shows $\pi_{\mathcal{H}} ({\n^t}_X Y) = \pi_{\mathcal{H}} ({\n}_X Y)$ and
$ \pi_\mathcal{V}(\n^t_XV)= \pi_\mathcal{V}({\n}_XV).$ The  first part of
the point (4) follows from the identities \eqref{A_and_Lie} and \eqref{A_altern} and the Koszul formulas. The last
equation follows from (1) and (2): $\n^t_UV={\n^t}^{\mathcal{V}}_UV +
T^t_UV={\n}^{\mathcal{V}}_UV + tT_UV= \n_UV+(t-1)T_UV.$
\end{proof}

\section{Twistor spaces over quaternionic and para-quaternionic K\"ahler manifolds}
In this section we consider pseudo-Riemannian submersions 
$\pi :(M,g) \ra (N,h) $ endowed with a complex structure $J$ on $M$  which is compatible with the decomposition \eqref{split_TM}. 
\bl\label{lemmaKaehler_subdata} Let $\pi :(M,g) \ra (N,h)$ be a
pseudo-Riemannian submersion endowed with a complex structure $J$ on $M$
which is compatible with the decomposition \eqref{split_TM}. Then $(M,g,J)$ is
a pseudo-K\"ahler manifold if and only if the following equations are satisfied 
\bea
&~&\pi_{\mathcal{H}} ((\n_XJ)Y) = \pi_{\mathcal{H}} ((\n_VJ)X)=0, \label{Kaehler_sub_data_I}\\
&~&(\n^{\mathcal{V}}_U	J)V=\pi_{\mathcal{V}} ((\n_XJ)V)=0,  \label{Kaehler_sub_data_II}\\
&~&A_X (JY) -JA_X Y =0, \quad A_X (JV) -JA_X V =0, \label{Kaehler_sub_data_III}\\
&~&T_V (JX) -JT_V X =0, \quad T_U (JV) -JT_U V=0,  \label{Kaehler_sub_data_IV}
\eea
where $X,Y$ are vector fields in $\mathcal H$ and $U,V$ are vector fields in $\mathcal V.$
\el
\begin{proof} Let $X,Y$ be vector fields in  $\mathcal H$ and $U,V$ be vector fields in $\mathcal V.$  Then it follows from comparing components in 
$\mathcal H\oplus \mathcal V $
\bean
(\n_XJ)Y&=& \pi_{\mathcal{H}} ((\n_XJ)Y) + \pi_{\mathcal{V}} (\n_XJ)Y) = \pi_{\mathcal{H}} ((\n_XJ)Y)+(A_X (JY) -JA_X Y),\\
(\n_XJ)V&=& \pi_{\mathcal{H}} ((\n_XJ)V) + \pi_{\mathcal{V}} (\n_XJ)V) = \pi_{\mathcal{V}} ((\n_XJ)V)+ (A_X (JV) -JA_X V),\\
(\n_VJ)X&=& \pi_{\mathcal{H}} ((\n_VJ)X) + \pi_{\mathcal{V}} (\n_VJ)X) = \pi_{\mathcal{H}} ((\n_VJ)X)+ (T_V (JX) -JT_V X),\\
(\n_UJ)V&=& \pi_{\mathcal{H}} ((\n_UJ)V) + \pi_{\mathcal{V}} (\n_UJ)V) =  (\n^{\mathcal{V}}_UJ)V+ (T_U (JV) -JT_U V).
\eean
\end{proof}
\noindent 
Further we define a second complex structure by
\bean \hat{J}:= 
\begin{cases} 
~J \mbox{ on } \mathcal{H},\\
-J \mbox{ on } \mathcal{V}.
\end{cases}
\eean
We observe that $\hat{\hat{J}}=J.$ This construction was made in \cite{N1} for the Riemannian setting and imitates the construction on twistor spaces.
\bp \label{twist_inv}
Suppose, that the foliation induced by the pseudo-Riemannian submersion $\pi $
is totally geodesic and that $(M,J,g)$ is a pseudo-K\"ahler manifold and $J$
is compatible with the decomposition \eqref{split_TM}, then the manifold
$(M,\hat{g}=g_{\frac{1}{2}},\hat{J})$ is a nearly pseudo-K\"ahler
manifold. The distributions $\mathcal{H}$ and $\mathcal{V}$ are parallel with
respect to the canonical Hermitian connection $\bar \n$ of $(M,\hat{g},\hat{J}).$
In other words the nearly pseudo-K\"ahler manifold $(M,\hat{g},\hat{J})$ has reducible holonomy. 
\ep
\pf Let $U,V$ be vector fields in $\mathcal{V}$ and $X,Y$ be vector fields in $\mathcal{H}:$  In the following $\hat \n$  is
the Levi-Civita connection of $\hat g.$ Since the fibers are totally geodesic,
i.e. $T\equiv 0,$ we obtain from equation \eqref{sub_zush_equ_I},
that  $\hat{\n}_UV = {\hat \n}^{\mathcal{V}}_UV +\hat T_UV={\n}^{\mathcal{V}}_UV +T_UV =\n_U V,$
  which yields $(\hat{\n}_U\hat{J})V= -({\n}_U{J})V=0.$ \\
In the sequel we denote the O'Neill tensors of the pseudo-Riemannian
foliations induced by $\mathcal V$ on $(M,{g})$ and on $(M,\hat{g})$ by
$A$ and $\hat{A},$ respectively.  From Lemma \ref{Besse_Lemma} it follows
$A_X Y = \hat{A}_X Y$ and consequently the same Lemma yields $\n_XY= \hat\n_XY.$ \\
Since $(M,g)$ is K\"ahler, Lemma \ref{lemmaKaehler_subdata} implies $A \circ J = J \circ A$ and we compute
\bea
(\hat{\n}_X\hat{J})Y&=& \hat{\n}_X(\hat{J}Y) -\hat{J}\hat{\n}_X Y\label{twist_inv_eqI}  \\ 
&=& \pi_{\mathcal{H}}[\hat{\n}_X({J}Y)] + \pi_{\mathcal{V}}[\hat{\n}_X({J}Y)]
-\hat{J}(\pi_{\mathcal{H}}(\hat{\n}_X Y) + \pi_{\mathcal{V}}(\hat{\n}_X
Y))  \nonumber\\ 
&=& \pi_{\mathcal{H}}[ \hat{\n}_X({J}Y) - {J}\hat{\n}_X Y ] + 
\pi_{\mathcal{V}}[\hat{\n}_X({J}Y) + J \hat{\n}_X Y] \nonumber \\
&=& \pi_{\mathcal{H}} ((\hat{\n}_XJ)Y) +\hat{A}_X({J}Y) + J \hat{A}_X Y \nonumber\\ 
&\overset{\eqref{Besse_Lemma_I_equ},\eqref{Besse_Lemma_III_equ},\eqref{Kaehler_sub_data_III}}=& \pi_{\mathcal{H}} (({\n}_XJ)Y) + 2 A_X(J Y)\overset{\eqref{Kaehler_sub_data_I}}=2 A_X(J Y)=2 J A_X Y.  \nonumber 
\eea
With the identity $A_X V = 2\hat{A}_X V$ of Lemma \ref{Besse_Lemma} we get
\bea
  (\hat{\n}_X\hat{J})V&=& \hat{\n}_X(\hat{J}V) -\hat{J}\hat{\n}_X V\label{twist_inv_eqII} \\
&=&  -\pi_{\mathcal{V}}(\hat{\n}_X({J}V)) -\pi_{\mathcal{H}}(\hat{\n}_X({J}V))+{J}\pi_{\mathcal{V}} (\hat{\n}_X V) - {J}\pi_{\mathcal{H}} (\hat{\n}_X V)\nonumber\\
&=&-\pi_{\mathcal{V}}((\hat{\n}_XJ)V)-  \hat{A}_X J V  - J \hat{A}_X V \nonumber \\
&\overset{\eqref{Besse_Lemma_I_equ},\eqref{Besse_Lemma_III_equ},\eqref{Kaehler_sub_data_III}}=& -\pi_{\mathcal{V}}(({\n}_XJ)V)-J {A}_X V = -{A}_X J V.\nonumber
\eea
The vanishing of the second fundamental form $T,$ equation \eqref{Besse_Lemma_Iprime_equ} and a second time 
 $A_X V = 2\hat{A}_X V$ show
\bea
(\hat{\n}_V\hat{J})X&=& 
\pi_\mathcal{V}({\hat \n}_V({J}X)) + \pi_\mathcal{V}({J}{\hat \n}_V X) 
+ \pi_\mathcal{H}(\hat{\n}_V({J}X) -{J}\hat{\n}_V X) \label{twist_inv_eqIII} \\\nonumber
&\overset{\eqref{Diff_n_nz_Hpart}}=& \hat T_V(J X) + J(\hat T_V X)+
\pi_{\mathcal H}((\n_VJ)X)+ \frac{1}{2}(J   A_X V-  A_{J X } V)=J A_X V,  
\eea
where we used $A_{J X}V=-J A_{X}V$ which follows, since $A_X$ is
alternating (compare equation \eqref{A_altern}) and commutes
with $J.$ The next Lemma finishes the proof.
\end{proof}

\bl \label{parallel_versus_sub_mersion_data}
1.) Suppose, that $(M,\hat J,\hat g)$ is a nearly pseudo-K\"ahler manifold and
$\hat J$ is compatible 
with the decomposition \eqref{split_TM}, then the following statements are equivalent:
\begin{itemize}
\item[(i)]
the splitting \eqref{split_TM} is $\bar \n$-parallel, 
\item[(ii)]
the fundamental tensors $\hat A$ and $\hat T$ satisfy: \\
\bea 
\hat T_VX&=&0, \quad \hat J \hat T_VW= -\hat T_V \hat J W \Leftrightarrow
\check J\hat T_VW= \hat T_V  \check J W \mbox{ for } \check J=\hat {\hat J},\\
\hat A_XV&=& \frac{1}{2}\hat J (\hat{\n}_X \hat J)  V, \quad 
\hat A_X Y =  \frac{1}{2}\pi_{\mathcal{V}}  \left( \hat{J} (\hat \n_X
\hat{J}) Y\right).
\eea
\end{itemize}
2.) If it holds $(\hat{\n}_V \hat J)  W=0$ then $\bar \n_VW \in \mathcal
V$ for $V,W \in \mathcal V$ is equivalent to $T_VW=0.$ Moreover it is
$(\hat{\n}^{\mathcal V}_V \hat J)  W=0.$ 
\el
\begin{proof}
First we compute 
\bean
\bar \n_VW &=& \hat{\n}_V W - \frac{1}{2}\hat J (\hat{\n}_V \hat J)  W \\
 &=& {\hat \n}^{\mathcal{V}}_V W + \hat T_V W - \frac{1}{2} \hat J ({\hat \n}^\mathcal{V}_V \hat J) W 
- \frac{1}{2}  (\hat J \hat T_V(\hat J W)+  \hat T_VW) \\
&=&  {\hat \n}^{\mathcal{V}}_V W  - \frac{1}{2} \hat J ({\hat \n}^\mathcal{V}_V \hat J) W + \frac{1}{2} (\hat T_VW-\hat J \hat T_V(\hat J W)) .
\eean
The first two terms lie in $\mathcal V,$  the second terms lie in $\mathcal{H}$ and therefore the expression is in
$\mathcal V$ if and only if  $\hat J\hat T_V( \hat J W)=  \hat T_VW
\Leftrightarrow \hat T_V(\hat  J W)= -\hat J
\hat T_VW.$ From   $0=\pi_{\mathcal H}(\bar \n_X V)= \pi_{\mathcal H}(\hat{\n}_X V - \frac{1}{2} \hat J (\hat{\n}_X \hat J)  V)$
one determines 
$$\hat A_XV=\pi_{\mathcal H}(\hat{\n}_X V) = \frac{1}{2} \pi_{\mathcal H}\left(\hat J (\hat{\n}_X \hat J)  V\right)=
\frac{1}{2} \hat J (\hat{\n}_X \hat J)  V.$$ 
The last equality follows from Lemma \ref{eig_V_T} (ii) part b). Conversely, if $\hat A_XV$ is given by the last formula one gets
\bean
\bar \n_X V=  \hat{\n}_X V - \frac{1}{2} \hat J (\hat{\n}_X \hat J)  V= \hat{\n}_X V -\hat A_XV=
  \pi_{\mathcal{V}}(\hat \n_X V) \in \mathcal{V}.
\eean
With the same identity we calculate
\bean
\pi_{\mathcal{V}}(\bar \n_V X) &=& \pi_{\mathcal{V}}\left(\hat{\n}_V X - \frac{1}{2} \hat J (\hat{\n}_V \hat J)  X\right) 
= \pi_{\mathcal{V}}(\hat{\n}_V X +\hat A_XV) = \pi_{\mathcal{V}} (\hat{\n}_V X).\eean
This is in $\mathcal H$ if and only if $\hat T_V X =\pi_{\mathcal{V}} (\hat{\n}_V X)=0.$ The last component, i.e. 
$\pi_{\mathcal{V}}(\bar \n_X Y) = \pi_{\mathcal{V}}\left(\hat \n_X Y - \frac{1}{2} \hat{J} (\hat \n_X  \hat{J}) Y\right)
$ is zero if and only if we have $\hat A_X Y = \pi_{\mathcal{V}}(\hat \n_X Y)
=\frac{1}{2} \pi_{\mathcal{V}}  \left( \hat{J} (\hat \n_X  \hat{J}) Y\right).$
\\
2.) With $(\hat{\n}_V \hat J)  W =0$ we calculate
\bean
\bar \n_VW &=& \hat{\n}_V W + \frac{1}{2} (\hat{\n}_V \hat J) \hat J W = {\hat
  \n}^{\mathcal{V}}_V W + \hat T_V W.
\eean
This lies in $\mathcal V$ if and only if $\hat T_VW=0.$
\end{proof}
\noindent
We apply Proposition \ref{twist_inv} to twistor spaces and obtain.
\bc \label{cor_tw_NK}
The twistor space $\mathcal{Z}$ of a quaternionic K\"ahler manifold of
dimension $4k$ with negative scalar curvature admits a canonical nearly pseudo-K\"ahler structure of
reducible holonomy contained in $U(1) \times U(2k).$
\ec
\pf
We remark that in negative scalar curvature the twistor space of a
quaternionic K\"ahler manifold is the total space of a pseudo-Riemannian
submersion with totally geodesic fibers. It admits a compatible
pseudo-K\"ahler structure of signature $(2,4k),$ cf. Besse \cite{Besse} 14.86
b). The assumption of positive scalar curvature is often made to obtain a
positive definite metric on $\mathcal{Z}.$ Here we focus on pseudo-Riemannian metrics 
and consequently on negative scalar curvature. 
\end{proof}

\bp 
The twistor spaces $\mathcal{Z}$ of non-compact duals of Wolf spaces and of 
Alekseevskian spaces admit a nearly pseudo-K\"ahler structure. 
\ep
\pf 
Non-compact duals of Wolf spaces are known \cite{W} to be quaternionic K\"ahler manifolds of
 negative scalar curvature. The same holds for Alekseevskian spaces \cite{A,C1}. \end{proof}
 
Studying the lists given in \cite{A,C1,W} we find examples of six-dimensional nearly pseudo-K\"ahler manifolds.
\bc
The twistor spaces $\mathcal{Z}$ of $\tilde {\mathbb{H}}P^1 = {\rm Sp}(1,1)/{\rm
  Sp}(1){\rm Sp}(1)$ and $SU(1,2)/S(U(1)U(2))$ provide six-dimensional nearly pseudo-K\"ahler manifolds.
\ec

\begin{remark} \label{REM_Example_twist}
The situation in negative scalar curvature is more flexible than in the positive
case. This is illustrated by the following results in this area: In the main
theorem of \cite{L} it is shown that the moduli space of complete
quaternionic K\"ahler metrics on $\mathbb{R}^{4n}$ is infinite dimensional. A
construction of super-string theory, called the {\it c-map} \cite{FS}, yields
continuous families of negatively curved quaternionic K\"ahler
manifolds. These results show that Corollary \ref{cor_tw_NK} is a good source
of examples.
\end{remark}

Another source of examples is given by twistor spaces over {\it para-quaternionic
K\"ahler manifolds}. Since these manifolds are less classical than quaternionic
K\"ahler manifolds, we recall some definitions (cf. \cite{AC} and references therein).

\bd \label{QKMFDEF}
Let $(\epsilon_1,\epsilon_2,\epsilon_3)=(-1,1,1)$ or some permutation thereof.
An {\cmssl almost para-quaternionic structure} on a differentiable manifold
$M^{4k}$ is a rank 3 sub-bundle $Q \subset End\, (TM),$ which is locally
generated by three anti-commuting endomorphism-fields $J_1,J_2,J_3=J_1J_2.$
These satisfy $J_i^2=\epsilon_i Id$ for $i=1,\ldots,3.$ Such a triple is
called {\cmssl   standard local basis} of $Q.$ 
A linear torsion-free connection preserving $Q$ is called {\cmssl
  para-quaternionic connection.} An almost para-quaternionic structure is
called a {\cmssl para-quaternionic   structure} if it admits a
para-quaternionic connection. An {\cmssl almost para-quaternionic hermitian
  structure}  $(M,Q,g)$ is a pseudo-Riemannian manifold endowed with a
para-quaternionic structure such that $Q$ consists of skew-symmetric
endomorphisms. For $n>1$ $(M^{4k},Q,g)$ is a {\cmssl para-quaternionic K\"ahler
  manifold} if $Q$ is preserved by the Levi-Civita connection of $g.$ In
dimension 4 a {\cmssl para-quaternionic K\"ahler  manifold} $M^4$ is an
anti-self-dual Einstein manifold.
\ed
We use the same notions omitting the word "para" for the quaternionic case. 
The condition that $Q$ is preserved by the Levi-Civita connection is
in a given standard local basis $\{J_i\}_{i=1}^3$ of $Q$  equivalent to
the equations 
\be \n_X J_i = -\theta_k(X)\epsilon_j J_j + \theta_j(X) \epsilon_kJ_k, \mbox{ for } X \in TM, 
\label{Q_parallel}\ee
where $i,j,k$ is a cyclic permutation of $1,2,3$ and $\{\theta_i \}_{i=1}^3$
are local one-forms. In the context of para-quaternionic manifolds one can
define twistor spaces for $s=1,0,-1$ 
$$ \mathcal Z^s:= \{ A \in Q \, | \, A^2=sId, \mbox{ with } A\ne 0\}.$$
The case of interest in this text is $\mathcal Z=\mathcal Z^{-1},$ since this twistor space is a complex manifold, such that the conditions
of Proposition \ref{twist_inv} hold true (cf. \cite{AC}). Therefore we obtain
the following examples of nearly pseudo-K\"ahler manifolds.
\bc The twistor space $\mathcal{Z}$ of a para-quaternionic K\"ahler
manifold with non-zero  scalar curvature of dimension $4k$  admits a canonical
nearly pseudo-K\"ahler structure of reducible holonomy contained in $U(k,k) \times U(1).$
\ec
\begin{example}
The para-quaternions $\widetilde{\mathbb{H}}$ are the $\mathbb{R}$-algebra
generated by $\{ 1, i, j, k\}$ subject to the relations $i^2=-1,~j^2=k^2=1,~
ij=-ji=k.$ Like the quaternions, the para-quaternions are a real Clifford
algebra which in the convention of \cite{LM} is
$\widetilde{\mathbb{H}}=\mathcal{C}l_{1,1}\cong\mathcal{C}l_{0,2}\cong\mathbb{R}(2).$
 One defines the para-quaternionic projective space $\widetilde{\mathbb{H}}P^n$ by the obvious equivalence
relation on the para-quaternionic right-module $\widetilde{\mathbb{H}}^{n+1}$
of $(n+1)$-tuples of para-quaternions.  The manifold
$\widetilde{\mathbb{H}}P^n$ is a para-quaternionic 
K\"ahler manifold \cite{Bl} in analogue to quaternionic projective space
$\mathbb{H}P^n.$ This yields examples of the type described in the last Corollary.
\end{example}
\section{Reducible nearly pseudo-K\"ahler manifolds}

In this section we study the case of a nearly pseudo-K\"ahler manifold
$(M^{2n},J,g),$ such that the holonomy of the canonical connection 
$\bar \n$ is reducible, in the sense that the tangent bundle $TM$ admits
 a splitting 
$$ TM = \mathcal{H} \oplus \mathcal{V}$$ 
into two $\bar \n$-parallel sub-bundles $\mathcal{H},\mathcal{V},$
 which are orthogonal and invariant with respect to the almost complex structure $J.$ 
\subsection{General properties}
In this subsection we carefully check, generalizing \cite{N2} to
pseudo-Riemannian foliations, the information which follows
from the decomposition into the $J$-invariant sub-bundles.
\bl \label{curv_lemma_dec} In the situation of this section and for a vector
field $X$ in $\mathcal H,$   a vector field $Y$ in $TM$ and vector fields 
$U,V$ in $\mathcal V$ it is
\bea \bar{R}(X,Y,U,V)=  g\left([ \n_U J,\n_V J ]X,Y \right)- g\left((\n_X J)Y,(\n_U J) V \right). \label{curv_id}\eea
\el
\begin{proof} 
Since $\mathcal{H}$ and $\mathcal{V}$ are $\bar \n$-parallel it follows
$\bar{R}(Y,U,X,V)=0$ and using equation \eqref{curv_rel_n_nbar} we get
\bean {R}(Y,U,X,V)&=&\frac{1}{2}g((\n_YJ)U,(\n_XJ)V) \\
&-& \frac{1}{4}\left[g((\n_YJ)X,(\n_{U}J)V)-  g((\n_YJ)V,(\n_{U}J)X)\right]\\
&=&-\frac{1}{2}g((\n_VJ)(\n_UJ)Y,X) \\
&-& \frac{1}{4}\left[g((\n_YJ)X,(\n_{U}J)V)-  g((\n_{U}J)(\n_VJ)Y,X)\right]. \eean
The first Bianchi identity yields 
\bean R(X,Y,U,V)&=&-R(Y,U,X,V) - R(U,X,Y,V)\\&=&-R(Y,U,X,V) + R(X,U,Y,V)\\
&=&\frac{3}{4}g\left([(\n_VJ),(\n_UJ)]Y,X\right) +\frac{1}{2} g((\n_YJ)X,(\n_{U}J)V).\eean
Replacing the last expression into 
\bean 
\bar{R}(X,Y,U,V)&=&  {R}(X,Y,U,V)-\frac{1}{2}g((\n_XJ)Y,(\n_UJ)V) \\
&+& \frac{1}{4}\left[g((\n_XJ)U,(\n_{Y}J)V)-  g((\n_XJ)V,(\n_{Y}J)U)\right]
\eean
proves the Lemma.
\end{proof}

\bc \label{cor_nablasqr}
For vector fields $X,Y$ in $\mathcal{H}$ and $V,W$ in $\mathcal{V}$ one has
\begin{enumerate}
\item[(i)] $ (\n_XJ)(\n_VJ)W=0;\, (\n_VJ)(\n_XJ)Y=0;$
\item[(ii)] $(\n_XJ)(\n_YJ)Z$ belongs to $\mathcal H$ for all $Z \in \Gamma(\mathcal{H});$
\item[(iii)] $(\n_VJ)(\n_WJ)X$ belongs to $\mathcal H;$ and  $(\n_XJ)(\n_YJ)V$ belongs to $\mathcal V.$ 
\end{enumerate}
\ec
\begin{proof}~\\
(i) follows from the fact, that $\bar{R}(J X,J Y,V,W)=\bar{R}(X,Y,V,W)$ and 
that the first term of equation \eqref{curv_id} has the same symmetry with
respect to $J.$ This yields on the one hand 
$$g\left((\n_{JX}J)JY,(\n_VJ) W \right)=g\left((\n_{X}J)Y,(\n_VJ) W \right)$$
and on the other hand it is 
 $$g\left((\n_{JX}J)JY,(\n_VJ) W \right)=-g\left((\n_{X}J)Y,(\n_VJ) W \right).$$
Consequently one has $g\left((\n_{X}J)Y,(\n_VJ) W \right)=0.$ Exchanging
$\mathcal H$ and $\mathcal V$ finishes part (i). \\
(ii) From (i) one gets the vanishing of
$$g((\n_VJ)(\n_YJ)Z,X)= g(Z,(\n_YJ)(\n_VJ)X)=-g(Z,(\n_YJ)(\n_XJ)V)=-g((\n_XJ)(\n_YJ)Z,V).$$
(iii) From (i) it follows $0=\bar{R}(X,U,V,W)=g([\n_VJ,\n_WJ]X,U).$ This yields 
$[\n_VJ,\n_WJ]X \in \mathcal{H}$ and by $[\n_VJ,\n_{JW}J]JX=-\{\n_VJ,\n_{W}J\}X\in \mathcal{H}$ we get the first part. The second part follows by replacing $\mathcal{H}$ and $\mathcal{V}.$
\end{proof}
\subsection{Co-dimension two}
Motivated by the above section on twistor spaces we suppose from now on that 
the real dimension of $\mathcal V$ is two.
\bl \label{eig_V_T} Let $\dim_{\bR}(\mathcal V)=2.$
\begin{itemize}
\item[(i)]
 Then the restriction of the metric $g$
  is either of signature $(2,0)$ or $(0,2).$ 
\item [(ii)]
\begin{itemize}
\item[a)]
 $T(V,W)=0$ for all $V,W \in \mathcal{V}.$
\item[b)]
 $T(X,U)\in \mathcal{H}$ for all $X \in \mathcal{H}$ and $U \in \mathcal{V}.$
\item[c)] In dimension six it is
 $T(X,Y)\in \mathcal{V}$ for all $X,Y \in \mathcal{H}.$ 
\item[d)]
 $\mbox{Span}\{\pi_{\mathcal{V}}(T(X,Y))\,|\, X,Y \in \mathcal{H}\}= \mathcal{V}.$ 
\end{itemize}
\end{itemize}
\el
\begin{proof}  Let $V \in \mathcal{V}$ with $g(V,V)\ne 0,$ then it is $\mathcal{V}= \mbox{span}\{ V, J V\}.$
(i) It holds $g(JV,JV)= g(V,V)\ne 0.$ This proves (i). \\
(ii) In the following we denote by $X,Y$ (local) sections of $\mathcal H$ and
by $U,V$ (local) sections of $\mathcal V:$ The part a) follows from $T(V,V)=T(JV,JV)=0$ and from the formula for the torsion $T(V,JV)= -J(\n_{V}J)JV=-(\n_{V}J)V=0.$ \\
Part b) follows from the fact that $g((\n_UJ)Y,V)$ and $g((\n_UJ)Y,JV)$ are three-forms:
This implies, that one has $g((\n_UJ)X,V)= -g((\n_UJ)V,X)\overset{a)}=0$ and 
$g((\n_UJ)X,JV)= -g((\n_UJ)V,JX)\overset{a)}=0.$ Hence the projection of $T(X,U)$ on $\mathcal{V}$ vanishes and part b) follows.
Next we show part d). There exists a pair $X,Y,$ such that
$\pi_{\mathcal{V}}(T(X,Y))\ne 0.$ Then it follows $\pi_{\mathcal{V}}(T(X,JY))=
-J\pi_{\mathcal{V}}(T(X,Y))\ne 0$ and part d) holds true. If there were no
such pair $X,Y$ then it follows $-g(T(X,Y),V)= g(J(\n_X J)Y,V)=-g(J(\n_V
J)Y,X)=0$ for all $X,Y \in \mathcal{H},$  $V \in \mathcal{V}$ and
consequently using a) and b) it follows $\n_VJ =0.$ This contradicts the fact, that $M$ is
strict nearly K\"ahler.  It remains part c). As $T(X,Y)\ne 0$ it follows that $X$ and $Y$ are
linear independent. From the symmetries of $T(X,Y)$ one concludes, that $\mbox{span}_{\bR}\{X,Y, JX,JY\}\subset \mathcal H$ has real dimension $4$ and
hence coincides with 
$\mathcal H.$ Using the symmetries of $\n J$ one gets
$T(X,Y),T(X,JY),T(JX,Y),T(JX,JY) \in \mathcal V.$ This finishes the proof.
\end{proof}
\bc 
Let $\dim_{\bR}(\mathcal V)=2.$ Then the foliation $\mathcal V$ has totally
geodesic fibers and the O'Neill tensor is given by $ A_XY = \frac{1}{2}
\pi_{\mathcal{V}}(J (\n_XJ)Y) $ and $ A_XV =\frac{1}{2} J (\n_XJ)V.$ Moreover
it is $ \n^{\mathcal V} J=0.$
\ec
\begin{proof}
From Lemma \ref{eig_V_T} (ii) a) we obtain $( \n_V {J})W=0$  with $V,W \in
\Gamma(\mathcal V).$ By Lemma \ref{parallel_versus_sub_mersion_data} part
2) it follows $T_VW=0$ and $ \n^{\mathcal V} {J}=0,$ since the decomposition $\mathcal{H} \oplus \mathcal{V}$ 
is  $\bar \n$ parallel.  Part 1) of Lemma \ref{parallel_versus_sub_mersion_data} finishes the proof. 
\end{proof}

\bp
Let $(M, J,  g)$ be a nearly pseudo-K\"ahler manifold  such that the property
of Lemma \ref{eig_V_T} (ii) c) is satisfied and such that $\mathcal{V}$ has dimension 2, then $(M, \check J=\hat{J},  \check g= g_2)$ is a pseudo-K\"ahler
manifold\footnote{Here we use $\check{\cdot}$ for the inverse construction of $\hat{\cdot}.$}.
\ep
It is natural to suppose the property of Lemma \ref{eig_V_T} (ii) c), since this 
holds true in the cases of {\it twistorial type} which are studied in the next
sections.
\begin{proof}
By the last Corollary the data of the submersion is $\check T=T\equiv 0,$ $
A_XY=\check A_XY = \frac{1}{2} \pi_{\mathcal{V}}( J (\n_X J)Y) $ and $
\check A_XV=2 A_XV =  J (\n_X J)V.$ Since $A$ anti-commutes
with $ J$ it commutes with $\check J.$ This yields the conditions
\eqref{Kaehler_sub_data_III} and \eqref{Kaehler_sub_data_IV} of Lemma
\ref{lemmaKaehler_subdata} on the triple $\check A,\check T, \check J.$ Further it holds $
\n^{\mathcal V}J=0.$ From the reasoning of equation \eqref{twist_inv_eqI} we
obtain $\pi_{\mathcal H}((\check\n_X\check J)Y)=\pi_{\mathcal H}((\n_X J)Y)$
which vanishes by the property of Lemma \ref{eig_V_T} (ii) c). By an analogous
argument we get from equation \eqref{twist_inv_eqII} the identity
$\pi_{\mathcal V}((\check\n_X\check J)V)=-\pi_{\mathcal V}((\n_X J)V).$ This
vanishes by Lemma \ref{eig_V_T} (ii) b). From equation \eqref{twist_inv_eqIII} we derive
$-\pi_\mathcal{H}((\n_XJ)V)\overset{n.K.}=\pi_\mathcal{H}((\n_VJ)X)=\pi_\mathcal{H}((\check \n_V \check J)X)
+2\pi_\mathcal{H}( J A_XV).$ The definition of $A_XV$ yields
$\pi_\mathcal{H}((\check \n_V\check J)X)=0.$
These are all the identities needed to apply Lemma \ref{lemmaKaehler_subdata}.	
\end{proof}

\bp \label{kurv_loc_sym}
Let $X,Y$ be vector fields in $\mathcal{H}$ and $V_1,V_2,V_3$  be vector fields
in $\mathcal{V}.$ 
Suppose that it holds $T(V,W)=0$ for all $V,W \in \mathcal{V}$ then it is 
\bea \label{barnablaformula}
\bar{R}((\n_XJ)JY,V_1,V_2,V_3)= g( JY,[\n_{V_1}J, [\n_{V_2}J,\n_{V_3}J]]X).
\eea
Moreover, one has $\bar\n_U \bar{R}(V_1,V_2,V_3,V_4)=0.$
\ep
\begin{proof} 
For $V_1,V_2,V_3 \in \mathcal V$ and $X\in \mathcal H$ the second Bianchi
identity gives
$$ -\underset{XYV_1}\sigma \bar \n_X (\bar R)(Y,V_1,V_2,V_3) =  \underset{XYV_1}\sigma \bar R((\n_XJ)JY,V_1,V_2,V_3).$$
As the decomposition $\mathcal H \oplus \mathcal V$ is $\bar \n$-parallel 
the terms on the left hand-side vanish due to the
symmetries \eqref{curiv_id_i} of the curvature tensor $\bar R.$ 
The right hand-side is determined  with the help of Lemma \ref{curv_lemma_dec} and Corollary \ref{cor_nablasqr}.
If we apply $\bar \n$ to the formula \eqref{barnablaformula} we obtain by
$\bar \n(\n J)=0$ the identity $ 
%g\left((\n_XJ) \left[\bar\n_U(\bar R)(V_1,V_2,V_3) \right],Y\right)=-
g\left(\bar\n_U(\bar R)(V_1,V_2,V_3) ,(\n_XJ)Z\right)=0$ with $Z=JY.$ This yields the proposition using Lemma \ref{eig_V_T} (ii) part d).
\end{proof}
\subsection{Six-dimensional nearly pseudo-K\"ahler manifolds}
Before analyzing the general case we first focus on dimension six.
\bl \label{leaf_curvature}
On a six-dimensional nearly pseudo-K\"ahler manifold $(M^6,J,g)$ the integral manifolds of the foliation $\mathcal{V}$ have Gaussian curvature $4\alpha$ 
and constant curvature $\kappa=4\alpha,$ where $\alpha$ is the type constant.
\el
\noindent
Recall that the sign of $\alpha$ is completely determined by the signature of
the metric $g,$ cf. section \ref{small_dim_sect}.
\pf
Let $X$ and $V$ be (local) vector fields of constant length in $\mathcal{H}$ and $\mathcal{V},$ respectively.
Then it follows from   equation \eqref{can_con_def} and the skew-symmetry  of $\n J$: $g(\n_VV, \n_X X)=g(\bar \n_VV,\bar \n_X X)=0.$
This identity yields
\bean
R(X,V,X,V)&=& -g(\n_V\n_X + \n_{[X,V]}V,X) \\
&\overset{NK}=&g(\n_XV,\n_VX) - \frac{1}{2}g(J(\n_{[X,V]}J)V,X)\\
&\overset{L. \ref{eig_V_T}\; (ii) \, b)}=&g(\n_XV,\n_VX) - \frac{1}{2}g(J(\n_XJ)V,\pi_{\mathcal{H}}([V,X]))\\
&\overset{(*)}=&g( \n_XV,\n_VX) -\frac{1}{2}g(J(\n_XJ)V,\n_VX-\frac{1}{2} J(\n_XJ)V)\\
&=& g(\bar \n_XV, \n_VX)+\frac{1}{4}g((\n_XJ)V,(\n_XJ)V)\\
&\overset{L. \ref{eig_V_T}\; (ii) \, b)}=& g(\bar \n_XV, \bar \n_VX)+\frac{1}{4}g((\n_XJ)V,(\n_XJ)V)\\
&=&\frac{1}{4}g((\n_XJ)V,(\n_XJ)V).
\eean
At the equality $(*)$ we use Lemma \ref{eig_V_T} (ii) b) which implies, that $\n_VX=\bar \n_VX  
+ \frac{1}{2} J(\n_VJ)X\in \mathcal{H}$ to show
\bean
 \pi_{\mathcal{H}}([V,X]) &=&\n_VX -\pi_{\mathcal H}(\bar \n_XV)- \frac{1}{2} J(\n_XJ)V \\
&=&\n_VX  - \frac{1}{2} J(\n_XJ)V.
\eean
Since $M^6$ is strict, we obtain $R(X,V,X,V)= \frac{1}{4}\alpha g(X,X) g(V,V).$
In addition it holds $Ric=5 \alpha g$ which implies 
\bean
 5\alpha g(V,V) &=&g(Ric(V),V) = g(V,V) R(JV,V,JV,V) +\sum_{i=1}^4 g(e_i,e_i) R(e_i,V,e_i,V)\\
&=& g(V,V) R(JV,V,JV,V) +\frac{\alpha}{4} \sum_{i=1}^4   g(e_i,e_i)^2 g(V,V),
\eean
where $\{e_i\}_{i=1}^4$ is an orthogonal frame of $\mathcal{H}.$
This equation yields $R(JV,V,JV,V)=4 \alpha$  and it follows that the fibers
have  Gaussian curvature
$$ K= \frac{R(JV,V,JV,V)}{g(V,V)g(JV,JV)-g(V,JV)^2}=4\alpha$$ and constant
curvature $\kappa= 4\alpha .$ 
\end{proof}

\bp
The manifold $(M,J,g)$ is the total space of a pseudo-Riemannian submersion
$\pi : (M,g) \ra (N,h)$ where $(N,h)$ is an almost pseudo-hermitian manifold
and the fibers are totally geodesic hermitian symmetric spaces. In particular,
the fibers are simply connected.
\ep
\pf
The foliation which is induced by $\mathcal{V}$ is totally geodesic and each
leaf is  by Proposition \ref{kurv_loc_sym} a locally hermitian symmetric space of complex dimension one.\\
It is shown in Lemma \ref{leaf_curvature} that each leaf has constant
curvature $\kappa.$ In the case $\kappa >0$ the leaves are compact and we can apply a
result of  Kobayashi, cf. \cite{Besse} 11.26, to obtain that the leaves are simply connected. Since the leaves
are also simply connected it follows, that the leaf holonomy is trivial and
that the foliation comes from a (smooth) submersion 
(cf. p. 90 of \cite{Sharpe}). In the case $\kappa < 0$ we observe, that  $(M,J,-g)$
is a nearly pseudo-K\"ahler manifold of constant type $-\alpha.$ The same
argument shows that the fibers are simply connected.
\end{proof}

\bl \label{Lem_quat_structure}
Let $(M^6, g,  J)$ be a strict nearly pseudo-K\"ahler six-manifold of constant
type $\alpha.$ For  an arbitrary normalized\footnote{Constant non-zero length
  suffices.} local vector field $V\in \mathcal{V},$ i.e. $\epsilon_V= g(V,V)\in \{\pm 1\},$ we
consider the endomorphisms $\tilde J_1:= J_{|\mathcal H},$ $\tilde{J}_2\,:\,
\mathcal{H} \ni X \mapsto (\n_V J)X \in \mathcal{H} $ and $\tilde J_3=\tilde J_1\tilde J_2.$
Then the triple $(\tilde J_1,\tilde J_2,\tilde J_3)$ defines an
$\epsilon$-quaternionic triple on $\mathcal H$  with $\epsilon_1=-1$ and $\epsilon_2=\epsilon_3=\mbox{sign}(-\alpha \epsilon_V)$ and it is \bean 
\pi^{\mathcal H}[(\n_\chi\tilde  J_i)Y]&=& - \theta_k(\chi)\,\epsilon_j \tilde
J_j Y +  \theta_j(\chi)\, \epsilon_k\tilde J_k Y,
\eean for a cyclic permutation of  i,j,k and with $\theta_1(\chi)=sign(\alpha)g(JV,\bar\n_\chi  V),$
$\theta_2(\chi)=-sign(\alpha)\sqrt{|\alpha|}g(V,J\chi)$ and
$\theta_3(\chi)=sign(\alpha)\sqrt{|\alpha|}g(V,\chi).$ 
The sub-bundle of endomorphisms spanned by  $(\tilde J_1,\tilde J_2,\tilde
J_3)$ does not depend on the choice of $V.$
\el
\begin{proof} 
Let $A(X)= (\n_V  J)X$ for a fixed $V\in \mathcal V$ with $\epsilon_V=g(V,V)\in \{\pm 1\}$ and an arbitrary $X \in \mathcal H.$
Then it is 
\bean  g(A^2(X),X)&=&- g(A(X),A(X))=- g((\n_V  J)X,(\n_V  J)X)\\
&=&-\alpha  g(V,V)  g(X,X)
\eean
and we find after polarizing the last expression in $X$ the identity $A^2(X) = -\alpha
\epsilon_V X.$ \\
Furthermore $A$ is a skew-symmetric endomorphism field and  in consequence
trace-free. Therefore the endomorphism field 
$$ \tilde J_2= 
\frac{1}{\sqrt{|\alpha|}}(\n_V  J) 
$$
is a hermitian structure if $\alpha \epsilon_V >0$  and a para-hermitian
structure (since it is trace-free) if $\alpha \epsilon_V <0.$ Next we set
$\tilde J_1= J_{| \mathcal H}$ and $\tilde J_3:=\tilde J_1 \circ \tilde J_2=-\tilde J_2\circ
\tilde J_1,$ which follows from $(\n J)\circ  J= - J\circ (\n  J)$ and observe ${\tilde J}_3^2={\tilde J}_2^2.$ Moreover these are (para-)hermitian structures, since
$ J$ and $\n_V J$ are skew-symmetric w.r.t. the metric $ g.$ Hence the
triple $(\tilde J_1,\tilde J_2,\tilde J_3)$ is a (para-)quaternionic triple on
$\mathcal H$ with $\epsilon_1=-1$ and $\epsilon_2=\epsilon_3=\mbox{sign}(-\alpha \epsilon_V).$ In the following we suppose, that  it holds
\bea \label{type_cond_pk_proof}
(\n_X \tilde J_1)Y \in \mathcal{V},\mbox{ for } X,Y \,\in \mathcal H. 
\eea
This identity yields $\pi^{\mathcal H}((\n_X\tilde J_1)Y)=0$
and in consequence it is 
\bean \pi^{\mathcal H}((\n_\chi\tilde J_1)Y)
&=& \epsilon_Vg(V,\chi) (\n_V\tilde J_1)Y+\epsilon_Vg(JV,\chi) (\n_{JV}\tilde
J_1)Y\\&=& \epsilon_Vg(V,\chi) (\n_V\tilde J_1)Y+ \epsilon_Vg(V,J\chi)
J(\n_{V}\tilde J_1)Y\\
&=& \epsilon_V\sqrt{|\alpha |}g(V,\chi) \tilde J_2 Y+
\epsilon_V\sqrt{|\alpha|}g(V,J\chi)\tilde J_3 Y \overset{!!}=-\theta_3(\chi) \epsilon_2\tilde J_2 Y+ \theta_2(\chi)\epsilon_3\tilde J_3 Y,
\eean
where we have to define
\bean \theta_2(\chi)&=&\epsilon_3\epsilon_V\sqrt{|\alpha|}g(V,J\chi)=-sign(\alpha)\sqrt{|\alpha|}g(V,J\chi), \\ 
\theta_3(\chi)&=&-\epsilon_2\epsilon_V\sqrt{|\alpha|}g(V,\chi)=sign(\alpha)\sqrt{|\alpha|}g(V,\chi).
\eean
Further we compute using the relation \eqref{can_con_def} for $\n$ and $\bar \n$ \be (\n_\chi  \tilde J_2) Y= \bar \n_\chi(  \tilde J_2 Y)-\tilde J_2 \bar
\n_\chi Y + \frac{1}{2\sqrt{|\alpha|}} \left[ J(\n_\chi J),(\n_V
   J) \right]Y, \mbox{ for } \chi \in TM, Y \in \mathcal
H \label{nabla_J_2_in_lemma} \ee
 and  get using $\bar \n (\n  J)=0$
\bean
\pi^{\mathcal H}[(\bar \n_\chi  \tilde J_2) Y] &=&\frac{1}{\sqrt{|\alpha|}} \pi^{\mathcal H}\left[(\bar \n_\chi (\n_V J)) Y\right] 
= - \frac{1}{\sqrt{|\alpha|}}\pi^{\mathcal H}[(\n_{\bar\n_\chi  V}  J)Y]\\
&=& -\frac{1}{\sqrt{|\alpha|}}\epsilon_V g(JV,\bar\n_\chi  V) (\n_{ JV}  J)Y = \frac{1}{\sqrt{|\alpha|}} \epsilon_V g(JV,\bar\n_\chi  V)J (\n_{V} J)Y,
\eean
where we recall that $\bar\n_\chi V \in \mathcal V$ has no part parallel
 to $V.$  Due to  equation \eqref{type_cond_pk_proof} and Lemma \ref{eig_V_T} (ii) the last
term of \eqref{nabla_J_2_in_lemma} lies in $\mathcal V$ if $\chi$ is in
$\mathcal H$ and vanishes if $\chi$ is a multiple of $JV.$ For $\chi=V$ we get
$\left[ J(\n_{V} J),(\n_V   J) \right]Y= 2|\alpha| \tilde J_3 \tilde
J_2Y=2 \epsilon_2 |\alpha|  \tilde J_1 Y.$ This shows
\bean \pi^{\mathcal H}[(\n_\chi  \tilde J_2) Y]&=&  \epsilon_2\sqrt{|\alpha|}
\epsilon_V g(V,\chi) \tilde J_1 Y + \epsilon_V g(JV,\bar\n_\chi  V) \tilde J_3 Y,\\
&=& - \theta_1(\chi)\epsilon_3 \tilde J_3 Y +  \theta_3(\chi)\epsilon_1 \tilde{J}_1Y,
\eean
if we set $\theta_1(\chi)=-\epsilon_3 \epsilon_V g(JV,\bar\n_\chi  V)=sign(\alpha)g(JV,\bar\n_\chi  V).$ It remains to differentiate the third (para-)complex structure:
\bea (\n_\chi  \tilde J_3) Y= \bar \n_\chi(  \tilde J_3 Y)-\tilde J_3 \bar
\n_\chi Y + \frac{1}{2\sqrt{|\alpha|}} [ J(\n_\chi  J),
   J(\n_V  J) ]Y. \label{nabla_J_3_in_lemma} \eea
Again one obtains using $\bar \n  J=0$ and $\bar \n (\n  J)=0:$ 
\bean
\pi^{\mathcal H}[(\bar \n_\chi  \tilde J_3) Y] &=& \frac{1}{\sqrt{|\alpha|}} \pi^{\mathcal H}[(\bar\n_\chi ( J \n_V J)) Y] = - \frac{1}{\sqrt{|\alpha|}} \pi^{\mathcal H}[ J(\n_{\bar\n_\chi
  V}  J)Y]\\
&=& - \frac{1}{\sqrt{|\alpha|}} \epsilon_V g(JV,\bar\n_\chi  V) J (\n_{ JV}
J)Y = \theta_1(\chi)\epsilon_2 \tilde{J}_2Y.
\eean
The last term of equation \eqref{nabla_J_3_in_lemma} lies (with the help of  equation \eqref{type_cond_pk_proof} and
Lemma \ref{eig_V_T} (ii)) in $\mathcal V$ for $\chi \in \mathcal H$ and
vanishes if $\chi$ is a multiple of $V.$ For $\chi=JV$ we compute
$\left[ J(\n_{JV} J),J(\n_V   J) \right]Y=\left[ (\n_{V} J),J(\n_V   J)
  \right]Y =  2|\alpha|\tilde J_2 \tilde J_3Y=-2\epsilon_3 |\alpha|  \tilde J_1
Y.$ This yields 
$$\pi^{\mathcal H}[(\n_\chi  \tilde J_3) Y]= - \theta_2(\chi)\epsilon_1 \tilde{J}_1Y + \theta_1(\chi)\epsilon_2 \tilde J_2 Y.$$
Given a second section $U$ in $\mathcal V$ with $g(U,U)= g(V,V)$ one has
$U=a V + b J V$ for real functions $a,b$ with $a^2 +b^2 =1.$ Using this one easily sees that
the triple induced by $V$ and the one by $U$ (locally) spans the same sub-bundle $Q$ of
endomorphisms of $\mathcal H.$ \end{proof}

\bl
Let $(M^6, g,  J)$ be a strict nearly pseudo-K\"ahler six-manifold of constant
type $\alpha.$ Let $s\,:\, U \subset N \ra M$ be a (local)
section\footnote{Local sections exist, since $\pi$ is locally trivial
  \cite{Besse} 9.3.} of $\pi$ on some open set $U.$ Define $\phi$ by  
$$\phi=s_*\circ \pi_* : \mathcal H_{s(n)} \overset{\pi_*}{\ra} T_nN
\overset{s_*}{\ra}s_*(T_nN) \subset T_{s(n)}M, \mbox{ for } n \in N$$
and set 
${J_i}_{|n}:= \pi_* \circ \tilde {J_i}_{|s(n)} \circ ({\pi_*}_{|\mathcal H})^{-1}$ for $i=1,\ldots,3,$
 where $\tilde J_i$ are defined in Lemma \ref{Lem_quat_structure}. Then
 $(J_1,J_2,J_3)$ defines a local $\epsilon$-quaternionic basis preserved by the
 Levi-Civita connection $\n^N$ of $N.$ 
\el
\begin{proof}
We choose $U$ such that the section $s$ is a diffeomorphism onto $W=s(U)$ and
a vector field  $V$ in $\mathcal V$ defined on a subset containing $W.$ As $\pi$ is a pseudo-Riemannian submersion we obtain from $\pi_* \circ
s_*=\id$ that $s$ is an isometry from $U$ onto $W.$ 
Therefore it holds $ s_*(\n^N_XY)=\pi^{s_*TN}[\n_{s_*X}s_*Y]$ which yields  $\n^N_XY=\pi_*(\n_{s_*X}s_*Y)$ 
and 
\bea ({\pi_*}_{|\mathcal H})^{-1}(\n^N_XY)&=& \pi^{\mathcal H}(\n_{s_*X}s_*Y).  \label{trafo_zush}
\eea 
For convenience let us identify $U$ and $W$ or in other words consider $s$ as
the inclusion $W\subset M.$ Then the projection on $s_*TN$ is
$\phi=s_*\pi_*=\pi_*|\mathcal{H}.$ Moreover we need the (tensorial) relation
$$ \n^N_X(\pi_* Z) - \pi_* \pi^{\mathcal H} (\n_X^MZ)=0 \mbox{ or equivalently }
\n^N_X \tilde Z -\pi_* \pi^{\mathcal H} (\n_X^M \phi^{-1} \tilde Z)=0, $$
which can be directly checked for basic vector fields. Using this identity  we get for $i=1,\ldots, 3$ 
\bean
\n^N_X(J_iY) &=& \n^N_X(\phi \, \tilde J_i \phi^{-1} Y) =  \phi \n^M_X( \tilde
J_i \phi^{-1} Y)= \phi \, (\n^M_X \tilde J_i) \, \phi^{-1} Y + \phi \,  \tilde
J_i \,  \n^M_X (\phi^{-1} Y)\\
&=& \phi \, (\n^M_X \tilde J_i) \, \phi^{-1} Y + \phi \,  \tilde
J_i \,  \phi^{-1} \n^N_X  Y=\phi \, (\n^M_X \tilde J_i) \, \phi^{-1} Y + J_i  \n^N_X  Y,
\eean
which reads $(\n^N_XJ_i)Y= \phi \, (\n^M_X \tilde J_i) \, \phi^{-1} Y.$ This
finishes the proof, since the right hand-side is completely determined by 
Lemma \ref{Lem_quat_structure}. Therefore we have checked the condition \eqref{Q_parallel}, i.e. the manifold $N$ is
endowed with a parallel skew-symmetric (para-)quaternionic structure, see also \cite{Besse} 10.32 and 14.36.
\end{proof} 
\subsection{General dimension}
In the last section we have seen that in dimension six the tensor $\n_V J$ 
induces a (para-)complex structure on $\mathcal{H}.$ This motivates the following
definition.
\bd \label{def_twist_type}
The foliation induced by $TM =\mathcal H \oplus \mathcal V$ is called of {\cmssl
  twistorial type} if for all $p \in M $ there exists a $V \in \mathcal V_p$ such
that the endomorphism 
$$ \n_V J \,:\, \mathcal H_p \ra \mathcal H_p$$
is injective.
\ed
\noindent
Obviously, if $\n_V J$ defines a (para-)complex structure, then the foliation
is of twistorial type.
\bp \label{prop_prop_A}
\begin{itemize}
\item[(a)]
If the metric induced on $\mathcal H$ is definite, then the foliation is of
twistorial type.
\item[(b)]
If the foliation is of twistorial type, then for all $p \in M $  and all $0\ne
U \in \mathcal V_p$ the endomorphism 
$$ \n_U J \,:\, \mathcal H_p \ra \mathcal H_p$$
is injective.
\item[(c)] It holds with $A:=\n_V J$ for some vector field $V$ in $\mathcal V$ of constant length and for vector fields $X \in \mathcal H$ and $\chi \in TM$
\be \bar \n_\chi(A^2)X=0 \label{der_Asquare}.\ee
Further it holds $[A^2,(\n_UJ)]=0$ for all $U \in \mathcal V$ 
and 
\be \n_U(A^2)X=0 \label{der_LC_Asquare}\ee 
for vector fields $U$ in $\mathcal V.$ 
\end{itemize}
\ep
\begin{proof}
Part (a) follows from $(\n_V  J) X \in \mathcal{H}$ for $X \in \mathcal H$ and
$V \in  \mathcal V,$ cf. Lemma \ref{eig_V_T} (i). For (b) we observe, that 
if $\n_V J$ is injective so is $\n_{ J V} J= - J \n_V  J.$ As $\mathcal V$ is of
dimension two $\{ V,  J V\}$ with $V \ne 0$ is an orthogonal basis. With $a,b \in \bR$
it follows $g((a\n_{V}J +b\n_{ J V} J)X,(a\n_{V} J +b\n_{ J
  V} J)X)
= (a^2+b^2)\, g((\n_{V} J)X,(\n_{V} J)X),$ which yields, that
$ \n_{aV+b  J V} J \,:\, \mathcal H_p \ra \mathcal H_p$
is injective since  $a\ne 0$ or $b\ne 0.$ It remains to prove part (c).
We first observe, that, since $V$ has constant length  and since $\bar
\n$ is a metric connection and preserves $\mathcal V,$ it follows $\bar
\n_\chi V = \alpha(\chi) J V$  for some one-form $\alpha.$ From
$\bar\n (\n J)=0$ we obtain
$$ (\bar \n_\chi A)X =(\bar \n_\chi(\n_VJ))X=(\n_{\bar \n_\chi V}J)X = \alpha(\chi)(\n_{J
  V}J)X =-\alpha(\chi)JAX$$
and  we compute using $\{A, J\}=0$
\bean  \bar \n_\chi(A^2)X= A(\bar \n_\chi A)X + (\bar \n_\chi A)AX
=-\alpha(\chi)[A(J(AX))+JA^2X]=0. \eean
The expression $[A^2,(\n_UJ)]=0$  is tensorial in $U$ and vanishes for $U=V.$
Therefore we only need to compute
$[A^2,(\n_{JV}J)]=-[A^2,J(\n_{V}J)]=-J[A^2,(\n_{V}J)] =0,$
where we used that $A^2$ commutes with $J.$ This implies 
\bean \n_U(A^2)X = \bar \n_U(A^2)X + \frac{1}{2}[J (\nabla_UJ), A^2]X=
-\frac{1}{2}[ (\nabla_{J U}J), A^2]X= 0 \eean
and proves part (c).
\end{proof}
\noindent In the following 
$V$ is a local vector field of constant length $\epsilon_V=g(V,V)\in \{\pm
1\}.$ \\
We denote by $\Omega$ the curvature form of the connection induced by $\bar
\n$ on the (complex) line bundle
$\mathcal V,$ which is given by 
$$ \bar{R}(X,Y)V= \Omega(X,Y)  J V, \mbox{ for } X,Y \in TM, V \in
  \mathcal V.$$
\bp \label{propo_asquare} If the foliation is of twistorial type, 
\begin{itemize}
\item[(i)]
then the endomorphism $A:=\n_V  J_{|\mathcal{H}}$ satisfies  $A^2= \kappa
\epsilon_V \id_{\mathcal H}$ for some real constant $\kappa\ne 0$
and $$\Omega= -2\kappa(2 \omega^{\mathcal V}-\omega^{\mathcal H}),$$ where
$\omega^{\mathcal H}(X,Y)= g(  X ,JY )$ is the restriction of
the fundamental two-form $\o$ to $\mathcal H;$
\item[(ii)]
for $X,Y$ in $\mathcal{H}$ it is  $(\n_X J)Y \in \mathcal V.$
\end{itemize}
\ep
\noindent
The proof of this proposition is divided in several steps. 
\bl
\begin{itemize}
\item[(i)]
For $X,Y$ in $\mathcal H$ and $V$ in $\mathcal V$ it is $\bar R(X,Y,V, J
V)=-2 g( (\n_V J)^2X, J Y).$
\item[(ii)]
For a given $X$ in $\mathcal H$ and $V$ in $\mathcal V$ it follows $\bar
R(X,V,V, JV)=0.$
\end{itemize}
\el
\begin{proof}
(i) Since $\mathcal H$ is $\bar \n$-parallel we obtain, that 
$\underset{XYV}\sigma \bar R(X,Y,V, J V)= \bar R(X,Y,V, J V).$
This is the left hand-side of the first Bianchi identity \eqref{curiv_id_iii} . The right hand-side
reads
\bean -\underset{XYV}{\sigma} g( (\n_X J)Y,
(\n_V J)  J V) &=& -g( (\n_V J)X,
(\n_Y J)  J V) -g( (\n_Y J)V,
(\n_X J)  J V)\\
&=&-2g( (\n_V J)^2X, J Y). 
\eean
\noindent
(ii) From the symmetries \eqref{curiv_id_i} of the curvature tensor $\bar R$ it follows 
$\bar R(X,V,V, JV)=\bar R(V, JV,X,V).$ This expression vanishes since $\mathcal
H$ is $\bar \n$-parallel.
\end{proof}
\noindent
From the last lemma we derive the more explicit expression of the curvature form $\Omega:$
\be \Omega = f \omega^{\mathcal V} +\epsilon_V\alpha, \label{Omega_expl} \ee
where $f$ is a smooth function, $\omega^{\mathcal V}$ is the restriction of the
fundamental two-form $\o=g(\cdot,J\cdot)$ to $\mathcal V$ and $\alpha(X,Y)=  - 2 g(A^2X,JY).$
\bl It holds with $U \in \mathcal V$ and $X,Y \in \mathcal H:$
\bea
 d\omega^{\mathcal V}(X,U,JU)&=&0 ,\label{curv_line_I_equ}\\
 d\alpha(X,U,JU)&=&0, \label{curv_line_II_equ}\\
d\omega^{\mathcal V}(U,X,Y)&=& -g(\n_UJ)X,Y),\label{curv_line_III_equ} \\
d\alpha(U,X,Y)&=&4 g(A^2(\n_UJ)X,Y).\label{curv_line_IV_equ}
\eea
\el
\begin{proof}
For vector fields $A,B,C $  on $M$ it is 
\bean
(\n_A\o^{\mathcal V})(B,C)
&=& A\o^{\mathcal V}(B,C)-\o^{\mathcal V}(\n_AB,C) -\o^{\mathcal V}(B,\n_AC)=
A\o^{\mathcal V}(B,C)\\ &~& -\o^{\mathcal V}\left(\bar
\n_AB+\frac{1}{2}J(\n_AJ)B,C\right)
-\o^{\mathcal V}\left(B,\bar \n_AC+\frac{1}{2}J(\n_AJ)C\right).
\eean
If two of them are  $X,Y \in \mathcal H$ and one is $U\in \mathcal V$ we check
using the definition of $\o^{\mathcal V},$ the information of Lemma
\ref{eig_V_T} and that the decomposition $\mathcal H \oplus \mathcal V$ is $\bar \n$-parallel:
\bean
&~&\n_U\o^{\mathcal V}(X,Y)=0, \\ 
&~&\n_X\o^{\mathcal V}(U,Y)=-\o^{\mathcal
  V}(U, \n_XY)= -\frac{1}{2}\o^{\mathcal
  V}(U, J(\n_XJ)Y),\\
&~&\n_Y\o^{\mathcal V}(X,U)=-\o^{\mathcal V}(\n_YX,U) =-\frac{1}{2}\o^{\mathcal V}(J(\n_YJ)X,U).
\eean
By the symmetries of $\o(J\n_\cdot J \cdot, \cdot)$ we conclude $d\o^{\mathcal
  V}(U,X,Y)=-g((\n_UJ)X,Y).$
Next we suppose $X \in  \mathcal H$ and $U \in \mathcal V$ and obtain with Lemma \ref{eig_V_T}:
\bean
&~&\n_X\o^{\mathcal V}(U,JU)= -\o^{\mathcal V}(\bar\n_X U, JU)- \o^{\mathcal
  V}(U,\bar \n_X(JU))=0, \\
&~&\n_U\o^{\mathcal V}(X,JU)=-\frac{1}{2}\o^{\mathcal V}(J(\n_U J) X,JU
)-\frac{1}{2}\o^{\mathcal V}(X,J(\n_U J) JU )=0,\\
&~&\n_{JU}\o^{\mathcal V}(U,X)=-\frac{1}{2}\o^{\mathcal V}(J(\n_{JU} J) U,X
)-\frac{1}{2}\o^{\mathcal V}(U,J(\n_{J U} J) X )=0.
\eean
This shows $d \o^{\mathcal V}(X,U,JU)=0.$ Let $X,Y \in \mathcal H$ and $U \in
\mathcal V.$ From $\alpha(U,\cdot)=0$ we conclude
\bean
(\n_U\alpha)(X,Y)&=&-2\left[g(\n_U(A^2)X,Y)+g(A^2X,(\n_UJ)Y)\right]\overset{\eqref{der_LC_Asquare}}=-2g(A^2X,(\n_UJ)Y),
\\
(\n_X\alpha)(U,Y)&=& -g(A^2(\n_UJ)X,Y),\;(\n_Y\alpha)(U,X)= -g(A^2(\n_UJ)Y,X),  
\eean
\noindent
which finishes the proof of $d\alpha(U,X,Y)=4g(A^2(\n_UJ)X,Y),$ since
$[A^2,\n_UJ]=0$ for all $U\in \mathcal V.$ We now prove the last identity
\bean
(\n_U\alpha)(JU,Y)&=& -\alpha(\bar \n_UJU +
\frac{1}{2}J(\n_UJ)JU,Y)-\alpha(JU,\bar \n_UY + \frac{1}{2}J(\n_UJ)Y)=0,\\
(\n_{JU}\alpha)(Y,U)&=& -\alpha(\bar \n_{JU}Y +
\frac{1}{2}J(\n_{JU}J)Y,U)-\alpha(Y,\bar \n_{JU}U +
\frac{1}{2}J(\n_{JU}J)U)=0,\\
(\n_{Y}\alpha)(U,JU)&=& -\alpha(\bar \n_{Y}U + \frac{1}{2}J(\n_{Y}J)U,JU)
-\alpha(U,\bar \n_{Y}JU + \frac{1}{2}J(\n_{Y}J)JU)=0,
\eean
where we used $\alpha(W,\cdot) =-\alpha(\cdot,W)= 0$ for $W \in \mathcal V.$  This finally shows  $d\alpha(X,U,JU)=0.$
\end{proof}

\begin{proof} (of the Proposition \ref{propo_asquare}) (i) Let $X,Y$ be vector fields in $\mathcal H$
  and $V$ be a local vector field in $\mathcal V$ of constant length. Since
  $\Omega$ as a curvature form of a (complex) line  bundle is closed, we
  obtain from equation \eqref{Omega_expl} 
$-\epsilon_Vd\alpha = fd\o^{\mathcal V} + df \wedge \o^{\mathcal V}.$
The equations \eqref{curv_line_I_equ} and \eqref{curv_line_II_equ} imply
$df_{|\mathcal H}=0.$ This implies $[X,Y]f=0$ and using that $\mathcal H$ 
is $\bar \n$-parallel we obtain $(\bar \n_X Y)f=0=(\bar \n_Y X)f$ which yields finally
$0= T^{\bar \n}(X,Y)(f)=-[ J(\n_X J)Y](f).$ By Lemma \ref{eig_V_T}
(ii) d) the last equation shows 
$df_{|\mathcal V}=0.$ Since $M$ is connected, it follows $f\equiv -\kappa$
for a constant $\kappa.$\\
Again using $d\Omega(V,X,Y)=0$   equation \eqref{curv_line_III_equ} and
 \eqref{curv_line_IV_equ} yield for arbitrary $X,Y$
$$ \kappa g((\n_V J)X,Y) +4 \epsilon_Vg(A^2 (\n_V J)X,Y)=0.$$
This implies $(\n_V J)(\kappa\id_{\mathcal H}+4\epsilon_V A^2)=0.$ Since the foliation is of twistorial type,
it follows $$A^2=-\epsilon_V \frac{\kappa}{4} \id_{\mathcal H}=-\epsilon_V \alpha
\id_{\mathcal H}$$ if we set $4\alpha= {\kappa}$ in analgogue to dimension six.\\
(ii) Since $\Omega$ is closed, it follows from part (i) and $d\o^{\mathcal
  V}(X,Y,Z)=0$ for $X,Y,Z \in  \mathcal H$ that it is $d \omega^{\mathcal H}(X,Y,Z) =0.$ Using
$d\omega^{\mathcal H}(X,Y,Z)= 3g((\n_X J)Y,Z)$ yields part (ii).
\end{proof}

\bp \label{Cor_quat_structure_dim_2n} 
Let $(M^{4k+2},g,J)$ be a strict nearly pseudo-K\"ahler manifold of
twistorial type. Let $s\,:\, U \subset N \ra M$ be a (local)
section of $\pi$ on some open set $U.$ Define $\phi$ by  
$$\phi=s_*\circ \pi_* : \mathcal H_{s(n)} \overset{\pi_*}{\ra} T_nN
\overset{s_*}{\ra}s_*(T_nN) \subset T_{s(n)}M, \mbox{ for } n \in N$$
and set ${J_i}_{|n}:= \pi_* \circ \tilde {J_i}_{|s(n)} \circ ({\pi_*}_{|\mathcal H})^{-1}$ for $i=1,\ldots,3,$ where $\tilde J_i$ are defined in Lemma \ref{Lem_quat_structure}. Then
 $(J_1,J_2,J_3)$ defines a local $\epsilon$-quaternionic basis preserved by the
 Levi-Civita connection $\n^N$ of $N.$ 
\ep
\begin{proof}
The proof of Proposition \ref{Lem_quat_structure} only uses $A^2= \kappa \epsilon_V  \id$
and $(\n_XJ)Y \in \mathcal{V}$ for $X,Y \in \mathcal{H}.$ Therefore we can generalize
it by means of Proposition \ref{propo_asquare} to strict nearly pseudo-K\"ahler manifolds of
twistorial type.
\end{proof}
\subsection{The twistor structure} \label{tw_proof}
In this subsection we finally characterize the nearly pseudo-K\"ahler
structures, which are related to the canonical nearly K\"ahler structure of
twistor spaces. 
\bt \label{NK_twist_last_thm}
\begin{enumerate}
\item[(i)]
The manifold $(M,J= \check J,\check g= g_2)$ is a twistor space of a quaternionic pseudo-K\"ahler
manifold, if it is $\epsilon_V \alpha  > 0.$
\item[(ii)]
The manifold $(M,J=\check J,\check g=  g_2)$ is a twistor space of a para-quaternionic K\"ahler
manifold, if it is $\epsilon_V \alpha  < 0.$
\end{enumerate}
\et
\begin{proof}
Denote by $\pi^{\mathcal Z} : \mathcal Z \ra N$ the twistor space of the
manifold $N$ endowed with the parallel skew-symmetric (para-)quaternionic
structure constructed from the foliation $\pi \,:\, M \ra N$ of twistorial type, cf. Proposition
\ref{Lem_quat_structure} for dimension six and Proposition \ref{Cor_quat_structure_dim_2n} for general dimension. We observe that the restriction of $ J$ to
$\mathcal H$ yields a (smooth) map
\bean 
\varphi \,:\, M \ra \mathcal Z, \quad m \mapsto  d\pi_m \circ
          J_m \,_{|\mathcal H} \circ 
(d{\pi_m}_{|\mathcal H})^{-1}=: j_{\pi(m)},
\eean
which by construction satisfies $\pi^{\mathcal Z}\circ  \varphi =  \pi$ and as a consequence
$d\pi^{\mathcal Z} \circ  d\varphi =  d\pi.$ Since $\pi$ and $\pi^{\mathcal Z}$
are pseudo-Riemannian submersions, the last equation implies that $d\varphi$
induces an isometry of the according horizontal distributions and maps the
vertical spaces into each other. Let us  determine the differential of $\varphi$ on $\mathcal V.$\\
{\bf Claim:}  For $V\in \mathcal V$ one has
\bean d\varphi(V) &=& 2  \; d\pi  \circ (\n_V J) \circ 
(d{\pi}_{|\mathcal H})^{-1},  \\  
d\varphi( JV) &=& 2 \;  d\pi \circ  (\n_{ J V} J) \circ 
(d{\pi}_{|\mathcal H})^{-1} =-2\; d\pi \circ   J(\n_V J) \circ 
(d{\pi}_{|\mathcal H})^{-1}. \eean
To prove the claim  we consider a (local) vector field $V\in \mathcal V$ and  a (local)
integral curve $\gamma$ of $V$ on some interval $I \ni 0$ with $\gamma(0)=m.$   Let $X$ be a vector field in $N.$
Denote by $\tilde X$ the horizontal lift of $X.$  The Lie transport of $\tilde X$ along the
vertical curve $\gamma$ projects to $X,$ i.e. it holds $d\pi_{\gamma(t)}(\tilde
X)=X$ for all $t\in I$ and in consequence $\left(d\pi_{\gamma(t)}\,_{|\mathcal
    H}\right)^{-1}X = \tilde X.$ In other words $d\pi$ commutes with this Lie
transport, which implies
\bean
d\varphi(V)X&=&  d\pi ( (\mathcal{L}_V J)\tilde X),
\eean
as one directly checks using basic vector fields. Therefore we need to determine the Lie-derivative $\mathcal{L}$ of $J:$
\bean
\pi^{\mathcal H}(({\mathcal L}_V  J)\tilde X) &=&
\pi^{\mathcal H}([V, J \tilde X]- J[V,\tilde X])\\
&=&\pi^{\mathcal H}\left( \n_V( J \tilde X)-\n_{ J \tilde X}V - J\n_V \tilde X
+ J \n_{\tilde X}V \right)\\
&=& \pi^{\mathcal H}\left( (\n_V J)\tilde X  -\frac{1}{2} J(\n_{ J \tilde X} J)V
+\frac{1}{2} J\left( J(\n_{\tilde X} J)\right)V \right)\\
&=&2 (\n_V J)\tilde X.
\eean
This shows $ d\varphi(V) = 2  \; d\pi  \circ (\n_V J) \circ (d{\pi}_{|\mathcal H})^{-1},$ which implies $ d\varphi( J V)= 2 \, d\pi  \circ (\n_{{ J} V} J)\circ (d{\pi}_{|\mathcal H})^{-1}=-2 \, d\pi  \circ J(\n_V J)\circ  (d{\pi}_{|\mathcal H})^{-1}.$ 
Given a local section $s \, :\, N \ra M$ and the associated adapted frame of
the (para-)quaternionic structure it follows that $\varphi \circ s$ is $J_1,$ $d\varphi(V)$ is
related to $J_2$ and  $d\varphi( J V )$ to $-J_3$ which span the tangent space of the fiber $F_{\pi(m)}=S^2$ in $\varphi(m).$
The complex structure of $\mathcal Z$ maps $J_2$ to $J_3.$ Hence $d\varphi$ 
is complex linear for the opposite complex structure $\check J$ on $M.$
Further one sees in this local frame that $\varphi$ maps horizontal part into horizontal part. Therefore $\varphi$ is an isometry
for the metric $\check g= g_2,$ i.e. the parameter $t=2$ in the canonical
variation of the metric $g.$ This means that $(M, \check J,\check g= g_2)$ is isometrically biholomorph to $\mathcal Z.$
\end{proof} 
\noindent
Combining Theorem \ref{split_dim_ten_thm} and Theorem \ref{NK_twist_last_thm}  we obtain {\bf Theorem A}.
%\bt 
%Let $(M^{10},J,g)$ be a nice decomposable nearly K\"ahler manifold,
%then the universal cover of $M$ is either the product of a pseudo-K\"ahler 
%surface and a (strict) nearly pseudo-K\"ahler manifold $M^6$ or it is a twistor
%space of an eight-dimensional (para-)quaternionic K\"ahler manifold endowed 
%with its canonical nearly pseudo-K\"ahler structure.
%\et

\end{document}